\documentclass{amsart}

\usepackage{amssymb, amsmath}
\usepackage{mathrsfs}
\usepackage{amscd}
\usepackage{verbatim}
\usepackage{enumerate}

\usepackage{hyperref}

\allowdisplaybreaks

\theoremstyle{plain}
\newtheorem{theorem}{Theorem}[section]
\theoremstyle{remark}
\newtheorem{remark}[theorem]{Remark}
\newtheorem{example}[theorem]{Example}
\theoremstyle{plain}
\newtheorem{corollary}[theorem]{Corollary}
\newtheorem{lemma}[theorem]{Lemma}
\newtheorem{proposition}[theorem]{Proposition}
\newtheorem{definition}[theorem]{Definition}

\numberwithin{equation}{section}

\def\Z{{\mathbb Z}}

\def\R{{\mathbb R}}
\def\C{{\mathbb C}}

\newcommand{\E}{{\mathbb E}}
\renewcommand{\P}{{\mathbb P}}
\newcommand{\F}{{\mathscr F}}

\renewcommand{\H}{{\mathscr H}}

\renewcommand{\a}{\alpha}

\newcommand{\g}{\gamma}
\renewcommand{\d}{\delta}
\newcommand{\e}{\varepsilon}
\renewcommand{\l}{\lambda}

\renewcommand{\O}{\Omega}

\newcommand{\calL}{{\mathscr L}}
\newcommand{\n}{\Vert}
\newcommand{\one}{{{\bf 1}}}
\newcommand{\embed}{\hookrightarrow}

\newcommand{\lb}{\langle}
\newcommand{\rb}{\rangle}

\newcommand{\limn}{\lim_{n\to\infty}}

\newcommand{\calA}{\mathscr{A}}

\newcommand{\nnn}{|\!|\!|}

\newcommand{\X}{\mathbb{X}^q}

\newcommand{\OO}{\mathcal{O}}
\renewcommand{\div}{{\rm div}\,}
\newcommand{\Xp}{X_{1-\frac1p,p}}
\newcommand{\XAp}{\Dom_A(1-\tfrac1p,p)}

\newcommand{\Dom}{\mathsf{D}}

\begin{document}

\author{Jan van Neerven}
\address{Delft Institute of Applied Mathematics\\
Delft University of Technology \\ P.O. Box 5031\\ 2600 GA Delft\\The
Netherlands} \email{J.M.A.M.vanNeerven@tudelft.nl}

\author{Mark Veraar}
\address{Delft Institute of Applied Mathematics\\
Delft University of Technology \\ P.O. Box 5031\\ 2600 GA Delft\\The
Netherlands} \email{M.C.Veraar@tudelft.nl}

\author{Lutz Weis}
\address{Institut f\"ur Analysis \\
Universit\"at Karlsruhe (TH)\\
D-76128  Karls\-ruhe\\Germany}
\email{Weis@math.uka.de}

\title
[Stochastic evolution equations]{Maximal
$L^p$-regularity for stochastic evolution equations}

\begin{abstract}
We prove maximal $L^p$-regularity for the stochastic evolution equation
$$
\left\{
\begin{aligned}
dU(t)  +  A U(t)\, dt& = F(t,U(t))\,dt + B(t,U(t))\,dW_H(t), \qquad t\in
[0,T],\\
 U(0) & = u_0,
\end{aligned}
\right.
$$
under the assumption that $A$ is a sectorial operator with
a bounded $H^\infty$-calculus
of angle less than $\frac12\pi$ on a space $L^q(\mathcal{O},\mu)$.
The driving process $W_H$ is a cylindrical Brownian
motion in an abstract Hilbert space $H$. For
$p\in (2,\infty)$ and $q\in [2,\infty)$
and initial conditions $u_0$
in the real interpolation space
$\XAp $
we prove existence of unique strong solution with trajectories in
$$ L^p(0,T;\Dom(A))\cap C([0,T];\XAp ),$$
provided the non-linearities
$F:[0,T]\times \Dom(A)\to L^q(\mathcal{O},\mu)$ and
$B:[0,T]\times \Dom(A) \to \g(H,\Dom(A^{\frac12}))$
are of linear growth and Lipschitz continuous in their second variables
with small enough
Lipschitz constants. Extensions to the case where $A$ is an adapted
operator-valued process are considered as well.

Various applications to stochastic partial differential
equations are worked out in detail. These include higher-order and time-dependent
parabolic equations and the
Navier-Stokes equation on a smooth bounded domain $\OO\subseteq \R^d$ with $d\ge 2$.
For the latter, the existence of a unique
strong local solution with values in $(H^{1,q}(\OO))^d$ is shown.
\end{abstract}

\subjclass[2000]{Primary: 60H15 Secondary: 35D10, 35R60, 46B09, 47D06, 47A60}

\keywords{Maximal $L^p$-regularity, stochastic evolution equations,
$R$-boundedness, $H^\infty$-functional calculus, stochastic Navier-Stokes equations}
\date\today

\thanks{The first named author is supported by VICI subsidy 639.033.604
of the Netherlands Organisation for Scientific Research (NWO). The second author
was supported by the Alexander von Humboldt foundation and VENI subsidy
639.031.930
of the Netherlands Organisation for Scientific Research (NWO).
The third named author is supported by a grant from the
Deutsche Forschungsgemeinschaft (We 2847/1-2).}

\maketitle

\section{Introduction}

Maximal $L^p$-regularity techniques have been pivotal in much of the recent progress
in the theory of parabolic evolution equations
(see \cite{Am, DHP, DoreVenni, KuWe, PrSo, We} and there references therein).
Among other things, such techniques provide a systematic and powerful tool to study
nonlinear and time-dependent parabolic problems.

For stochastic parabolic evolution equations, maximal $L^p$-regularity
results have been obtained previously by Krylov
for second order problems on $\R^d$ \cite{Kry94a, Kry96, Kry, Kry00, Kry06},
by Kim for second order problems on bounded domains in $\R^d$ \cite{Kim09},
and by Mikulevicius and Rozovskii for Navier-Stokes equations \cite{MiRo04}.
A systematic theory of maximal $L^p$-regularity for stochastic evolution equations, however,
based on abstract operator-theoretic properties of the operators governing the equation,
has yet to be developed.
A first step towards such a theory has been taken in our recent paper \cite{NVW10},
where it was shown that if $A$ is a sectorial operator with a bounded $H^\infty$-calculus of angle $<\frac12\pi$
on a space $L^q(\mathcal{O},\mu)$ with $(\mathcal{O},\mu)$
an arbitrary $\sigma$-finite measure space and $q\in [2,\infty)$,
then $A$ has stochastic maximal $L^p$-regularity for all
$p\in (2,\infty)$, i.e., $A$ satisfies the convolution
 estimate
\begin{equation}\label{eq:maxregineq}
 \Big\n t\mapsto \int_0^t A^\frac12 S(t-s)G(s)\,dW_H(s) \Big\n_{L^p(\R_+\times\O;L^q(\mathcal{O},\mu))}
\le C \n G\n _{L^p(\R_+\times\O;L^q(\mathcal{O},\mu;H))},
\end{equation}
where $S$ denotes the semigroup generated by $-A$
and $W_H$ is a cylindrical Brownian
motion in a Hilbert space $H$. The stochastic integral is
understood as a vector-valued stochastic integral in $L^q(\mathcal{O},\mu)$
in the sense of
\cite{NVW1}.

The aim of this paper is to apply the above estimate to deduce maximal $L^p$-regularity
for the stochastic parabolic evolution equation
$$
\left\{
\begin{aligned}
dU(t)  +  A U(t)\, dt& = F(t,U(t))\,dt + B(t,U(t))\,dW_H(t), \qquad t\in
[0,T],\\
 U(0) & = u_0,
\end{aligned}
\right.
$$
Our main result asserts that if $A$ has a bounded $H^\infty$-calculus
of angle $<\frac12\pi$ on a Banach space $X$ that is isomorphic to a
closed subspace of $L^q(\mathcal{O},\mu)$ with $q\in [2,\infty)$, then for
$p\in (2,\infty)$ and initial conditions $u_0$
in the real interpolation space $\XAp  = (X,\Dom(A))_{1-\frac1p,p}$,
this problem
has a unique strong solution with trajectories in
$$ L^p(0,T;\Dom(A))\cap C([0,T];\XAp ),$$
provided the non-linearities
$F:[0,T]\times \Dom(A)\to X$ and
$B:[0,T]\times \Dom(A) \to \g(H,\Dom(A^{\frac12}))$
are of linear growth and Lipschitz continuous in their second variables with small enough
Lipschitz constants. The precise statement is contained in Theorem \ref{thm:SE},
where we allow $A$, $F$ and $B$, $u_0$ to be random.

To illustrate the power of this result, we apply it to the time-dependent
problem
$$
\left\{
\begin{aligned}
dU(t)  +  A(t) U(t)\, dt& = F(t,U(t))\,dt + B(t,U(t))\,dW_H(t), \qquad t\in
[0,T],\\
 U(0) & = u_0,
\end{aligned}
\right.
$$
and show in Theorem \ref{thm:SE2} that, essentially
under the same assumptions as in the time-independent case,
the same conclusions can be drawn with regard to the existence, uniqueness, and
regularity of strong solutions. An extension to the case of locally Lipschitz continuous
coefficients is given in Subsection \ref{subsec:local}.
These results extend \cite[Theorems 4.3 and 4.10]{Brz2}, \cite[Theorem 2.5]{Zha} and
\cite[Theorem 6.1]{Zha10} to the
case of sharp exponents.

It has already been mentioned that in Theorem \ref{thm:SE} we allow $A$ to be random.
In the special case where
$A$ is a fixed deterministic operator, the theorem can be applied (by taking
the negative extrapolation space $\Dom(A^{-\frac12})$ as the state space) to the situation
where the non-linearities are of the form
$F:[0,T]\times \Dom(A^\frac12)\to \Dom(A^{-\frac12})$ and
$B:[0,T]\times \Dom(A^\frac12) \to \gamma(H,X)$.
For initial values in $\Dom_A(\frac12-\frac1p,p)$, this results in solutions
with trajectories in  $L^p(0,T;\Dom(A^\frac12))
\cap C([0,T]; \Dom_A(\frac12-\frac1p,p)).$
For second order elliptic operators $A$ on a smooth domain  $\OO\subseteq \R^d$,
this includes the case where $F$ and $B$ arise as Nemytskii operators
 associated with nonlinear functions of the form $f(u,\nabla u)$ and $b(u,\nabla(u))$.
This is because in this setting $\Dom(A^\frac12)$ typically can be identified as a
Sobolev space $H^{1,q}$.
An illustration is given in Section \ref{sec:NS}, where we prove existence of
solutions in $H^{1,q}$ for the stochastic Navier-Stokes equation.

The advantage of the abstract approach presented in this paper
is that it replaces some of the hard (S)PDE techniques
of Krylov's $L^p$-theory by the generic assumption that $A$ have a good functional calculus.
In recent years, a large body of results has been accumulated by many authors
which shows that, as a rule of thumb, any `reasonable' elliptic operator of order $2m$
has such a calculus (see \cite{AHS, DDHPV, DHP, DuMc, DuRob, DuSi, DuYa, KKW, KWcalc, KuWe, McI, NoSa, Weis-survey}
 and the references therein); much of the hard analysis goes into proving these ready-to-use results.
Moreover, in most of these examples,
the trace space $\XAp$ and the fractional domain space $\Dom(A^\frac12)$
have been characterised explicitly as a fractional Besov space of order
$2m(1-\frac1p)$ and a Sobolev space of order $m$, respectively.

\subsection{Applications}
In principle, our results pave the way for proving maximal $L^p$-regularity
results for any parabolic problem governed by an operator having a
bounded $H^\infty$-calculus.

To keep this paper at a reasonable length we have picked
three examples which we believe to be representative (but by no means exhaustive)
to illustrate the scope of applications.
Further potential applications include, for instance,
parabolic SPDEs on complete Riemannian manifolds and on Wiener spaces such as
considered in \cite{Zhang07}
(cf. Examples \ref{ex:Hinfty} (7) and (8) below).

\subsubsection{Higher-order parabolic SPDEs on $\R^d$}

Our first application concerns a system of $N$ coupled parabolic SPDEs
involving elliptic operators of order $2m$ on $\R^d$ of the form
\begin{equation*}
\left\{
\begin{aligned}
d u(t,x) +  \mathcal{A}(t,x, D) u(t,x)\, dt &= f(t,x,u)\,dt +
\sum_{i\geq 1} b_i(t,x,u) \, d w_i(t),
\\ u(0,x) & = u_0(x).
\end{aligned}
\right.
\end{equation*}
Here
\[\mathcal{A}(t,\omega,x, D) = \sum_{|\alpha|\leq 2m} a_{\alpha}(t,\omega,x)
D^{\alpha}\]
with $D = -i(\partial_1, \ldots, \partial_d)$.
The scalar Brownian motions $w_i$ are independent,
and the functions $f$ and $b_i$
are Lipschitz continuous with respect to the graph norm of
$\mathcal{A}$.
Under suitable boundedness and continuity assumptions on the coefficients
$a_{\a}$ and a smallness condition on Lipschitz constants of $f$ and $b_i$
we prove the existence and uniqueness of a strong solution
with values in $H^{2m,q}(\R^d;\C^N)))$ and with continuous trajectories
in the Besov space $B^{2m(1-\frac1p)}_{q,p}(\R^d;\C^N))$ (Theorem \ref{thm:eq2m}).
To the best of our knowledge, this is the first maximal $L^p$-regularity result
for this class of equations.

\subsubsection{Time-dependent second-order parabolic SPDEs on bounded domains}
As a second example we consider time-dependent parabolic second order problems
on a bounded domain $\OO\subseteq \R^d$ whose boundary consists of two disjoint arcs
$\partial \OO = \Gamma_0\cup \Gamma_1$. We impose Dirichlet
conditions on $\Gamma_0$ and Neumann conditions on $\Gamma_1$
and prove the existence of a unique strong solution with values
in $H^{2,q}(\mathcal{O})$
and with continuous trajectories
in the Besov space $B^{2-\frac2p}_{q,p}(\mathcal{O})$
(Theorem \ref{thm:eq2mneumann}).

\subsubsection{The Navier-Stokes equation on bounded domains}
In the final section
we consider the stochastic Navier-Stokes equation in a bounded smooth domain $\OO\subseteq\R^d$ with $d\ge 2$
subject to Dirichlet boundary conditions. We prove existence and uniqueness
of a local mild solution with values
in $(H^{1,q}(\OO))^d$ and with continuous trajectories in $(B^{1-\frac2p}_{q,p}(\OO))^d$
for $\frac{d}{2q} < 1-\frac2p$.

\section{Preliminaries}

The aim of this section is to fix notations and to recall some recent results
on maximal $L^p$-regularity and stochastic maximal $L^p$-regularity that will
be
needed in the sequel.

Throughout this article we fix a probability spaces $(\O,\calA,\P)$
endowed with filtration $\F = (\F_t)_{t\ge 0}$, a Hilbert space $H$
with inner product $[\cdot,\cdot]$, and a Banach space $X$.

For $p_1, p_2\in [1,\infty]$,
the closed linear span
in $L^{p_1}(\O;L^{p_2}(\R_+;X))$ of all processes of the form
$f = \one_{(s,t]\times F}\otimes x$ with $F\in \F_s$ and $x\in X$
is denoted by
$$L^{p_1}_{\F}(\O;L^{p_2}(\R_+;X)).$$
The elements in $L_{\F}^{p_1}(\O;L^{p_2}(\R_+;X))$ will be referred to as the {\em
$\F$-adapted}
elements in $L^{p_1}(\O;L^{p_2}(\R_+;X))$.

The vector space of all (equivalence classes of)
strongly measurable functions on $\O$ with values
in a Banach space $Y$ is denoted by $L^0(\O;X)$.
The topology of convergence in probability is metrised by
the distance function $d(f,g) = \E (\n f - g\n\wedge 1)$
which turns $L^0(\O;Y)$ into a
complete metric vector space. The space of all $f\in L^0(\O;Y)$
that are strongly $\mathscr{B}$-measurable, where $\mathscr{B}\subseteq\calA$ is a
sub-$\sigma$-algebra, is denoted by $L_{\mathscr{B}}^0(\O;Y)$.

\subsection{Stochastic integration}
We will be interested in an estimate for stochastic integrals of the form
$\int_{\R_+} G\,dW_H,$
where $G$ is an $\F$-adapted process
with values in space of finite rank operators from $H$ to $X$,
and $W_H$ is an $\F$-cylindrical Brownian motion in $H$. We start with a concise
explanation of these notions.

\subsubsection{The space $\gamma(\H,X)$}
Let $\H$ be a Hilbert spaces (typically we take $\H=H$ or $\H = L^2(\R_+;H)$).
The space of all $\gamma$-radonifying operators from $\H$ to $X$ is denoted by
$\gamma(\H,X)$. Recall that this space is the closure of the space of
finite rank operators from $\H$ to $X$ with respect to the norm
$$ \Big\n \sum_{n=1}^N h_n\otimes x_n\Big\n_{\gamma(\H,X)}^2:=
\E  \Big\n \sum_{n=1}^N \gamma_n\otimes x_n\Big\n^2,$$
where it is assumed that $(h_n)_{n=1}^N$ is an orthonormal sequence in $\H$,
$(x_n)_{n=1}^N$ is a sequence
in $X$, and $(\gamma_n)_{n=1}^N$ is any sequence of independent standard Gaussian random variables.
For expositions
of the
theory of $\gamma$-radonifying operators we refer to \cite{DJT} and the review article \cite{NeeCMA},
where also references to
the extensive literature can be found.

For $X = L^p(\mathcal{O},\mu)$ with $1\le p,\infty$ and $(\mathcal{O},\mu)$ $\sigma$-finite,
one has a canonical isomorphism
\begin{align}\label{eq:gammaFubini0}
 L^p(\mathcal{O},\mu;\H) \simeq \gamma(\H, L^p(\mathcal{O},\mu))
\end{align}
which is obtained by assigning to a function $f\in L^p(\mathcal{O},\mu;\H)$
the operator $T_f: H\to L^p(\mathcal{O},\mu)$, $h\mapsto [f(\cdot),h]$ (see \cite{BrzvN03}).
More generally the same procedure gives, for any Banach space $X$,
a canonical isomorphism
\begin{align}\label{eq:gammaFubini}
L^p(\mathcal{O},\mu;\gamma(\H,X)) \simeq \gamma(\H, L^p(\mathcal{O},\mu;X))
\end{align}
(see \cite{NVW1}). We shall need the following variation on this theme.
Recalling the definition of the Bessel potential
spaces $H^{2\alpha,p}(\mathcal{O})$, where $\mathcal{O}\subseteq \R^n$ is a smooth domain,
application of the operator $(I-\Delta)^{-\a}$
on both sides of \eqref{eq:gammaFubini} gives an isomorphism
\begin{align}\label{eq:BesselFubini}  H^{2\a p}(\mathcal{O};\gamma(\H,X))
\simeq \gamma(\H, H^{2\a,p}(\mathcal{O};X)).
\end{align}

\subsubsection{Cylindrical Brownian motions}
An {\em $\F$-cylindrical Brownian motion in $H$} is a bounded linear
operator $W_H: L^2(\R_+;H)\to L^2(\O)$ such that:
\begin{enumerate}[\rm(i)]
\item for all $f\in L^2(\R_+;H)$ the random variable
$W_H(f)$ is centred Gaussian.
\item for all $t\in \R_+$ and $f\in L^2(\R_+;H)$ with support in $[0, t]$, $W_H(f)$ is $\F_t$-measurable.
\item for all $t\in \R_+$ and $f\in L^2(\R_+;H)$ with support in $[t, \infty)$, $W_H(f)$ is independent of $\F_t$.
\item for all $f_1,f_2\in L^2(\R_+;H)$ we have
$ \E (W_H(f_1)\cdot W_H(f_2)) = [f_1,f_2]_{L^2(\R_+;H)}.$
\end{enumerate}
It is easy to see that for all $h\in H$ the process $(W_H(t)h )_{t\ge 0}$ defined by
$$W_H(t)h := W_H(\one_{(0,t]}\otimes h)$$
is an $\F$-Brownian motion (which is standard if $\n h\n=1$). Moreover, two such
Brownian motions
$((W_H(t)h_1)_{t\ge 0}$  and $((W_H(t)h_2)_{t\ge 0}$ are independent if and only
if $h_1$ and $h_2$ are orthogonal in $H$.

\begin{example}[Space-time white noise]
Any space-time white noise $W$ on a domain $\OO\subseteq \R^d$ defines
a cylindrical Brownian motion in $L^2(\OO)$ and vice versa by the formula
$$W_{L^2(\OO)}(\one_{(0,t]}\otimes \one_B) = W(t,B)$$
for Borel sets $B\subseteq \OO$ of finite measure.
\end{example}

\begin{example}[Sums of independent Brownian motions]
A family $(w_i)_{i\in I}$ of independent
real-valued standard Brownian motions defines a cylindrical Brownian motion
in $\ell^2(I)$ and vice versa by
$$ W_{\ell^2(I)}(\one_{(0,t]}\otimes e_i) := w_i(t),$$
where $e_i\in \ell^2(I)$ is given by $e_i(j) = \delta_{ij}$.
\end{example}

\subsubsection{The stochastic integral}
Processes which are finite linear combinations of processes of the form
$$\one_{(s,t]\times F}\otimes (h\otimes x)$$ with $F\in \F_s$, $h\in H$, $x\in
X$,
are called {\em $\F$-adapted finite rank step processes in $\gamma(H,X)$}.
The stochastic integral of such a process
with respect to an $\F$-cylindrical Brownian motion $W_H$ is defined
 by
$$  \int_{\R_+} \one_{(0,t]\times F}\otimes (h\otimes x) \,dW_H:=
\one_F [W_H(t)h]\otimes x $$
and linearity.
The following two-sided estimate has been proved in \cite{NVW1}:

\begin{theorem}\label{thm:UMD}
Let $X$ be a UMD Banach space and let $G$ be an $\F$-adapted
finite rank step process in $\gamma(H,X)$. For all $p\in (1,\infty)$
one has the two-sided estimate
\begin{align}\label{eq:twosided} \E\Big\|\int_{\R_+} G(s) \, dW_H(s)\Big\|^p \eqsim_p
\E\|G\|_{\gamma(L^2(\R_+;H),X))}^p,
\end{align}
with implicit constants depending only on $p$ and (the UMD constant of) $X$.
\end{theorem}

This equivalence is used to give a meaning to the stochastic integral on the
left-hand side of the maximal $L^p$-regularity inequality \eqref{eq:maxregineq}
and plays a crucial role in the proof of this inequality; the inequality \eqref{eq:type2est}
does not suffice for this purpose (see \cite{NVW10}).

Examples of UMD spaces are all Hilbert spaces and the
spaces $L^q(\OO,\mu)$ with $q\in (1, \infty)$. Furthermore,
closed subspaces, quotients, and duals of
UMD spaces are UMD.  For more information on UMD spaces we
refer to \cite{Bu3}.

As a consequence of Theorem \ref{thm:UMD} and a routine density argument,
the stochastic integral can be uniquely extended to the space
$L_\F^p(\Omega;\gamma(L^2(\R_+;H),X))$, which
is defined as the closed linear span in $L^p(\Omega;\gamma(L^2(\R_+;H),X))$ of all
$\F$-adapted finite rank step processes in $\gamma(H,X)$. For a detailed discussion
we refer to \cite{NVW1}.

For Banach space $X$ with type $2$ one has a continuous embedding
\begin{align*}
L^2(\R_+;\g(H,X))\embed \gamma(L^2(\R_+;H),X)
\end{align*}
(see \cite{NW1, RS}). In combination with \eqref{eq:twosided} this gives the following estimate,
valid for finite rank step process in $\gamma(H,X)$ with $X$ a UMD space with type $2$:
\begin{align}\label{eq:type2est}
\E\Big\|\int_{\R_+} G(s) \, dW_H(s)\Big\|^p \leq C^p
\E\|G\|_{L^2(\R_+;\g(H,X))}^p.
\end{align}
As a consequence of the inequality \eqref{eq:type2est}, the stochastic integral uniquely extends to
$L_\F^p(\O;L^2(\R_+;\g(H,X)))$,
the closed linear span in $L^p(\O;L^2(\R_+;\g(H,X)))$ of all
$\F$-adapted finite rank step processes in $\gamma(H,X)$.

Examples of UMD spaces with type $2$ are all Hilbert spaces and the
spaces $L^q(\OO,\mu)$ with $q\in [2, \infty)$. A
UMD space has type $2$ if and only if it has martingale type $2$,
and in fact the estimate \eqref{eq:type2est} holds for any Banach space $X$ with
martingale type $2$ (see \cite{Brz2, Nh}).
For more information on the notions of (martingale) type and cotype we
refer to \cite{DJT, Pi75, Pi2}.

\begin{remark}
It follows easily from \cite{Lenglart} that the estimates \eqref{eq:twosided} and \eqref{eq:type2est} are valid for
arbitrary exponents $p\in (0,\infty)$.
We shall not need this fact here.
\end{remark}

\subsubsection{The stochastic integral operator family $\mathscr{J}$}\label{sec:J}
We turn our attention to a class of stochastic integral operators, which plays a key role in connection with
stochastic maximal $L^p$-regularity (see Theorem \ref{thm:maxregstoch} below).

For an $\F$-adapted finite rank step process
$G:\R_+\times\O \to \g(H,X)$
and a parameter $\delta>0$ we define the process
$J(\delta)G:\R_+\times\O \to X$ by
\[(J(\delta)G)(t) := \frac1{\sqrt{\d}} \int_{(t-\delta)\vee 0}^t  G(s) \,
dW_H(s).\]

A routine computation using
\eqref{eq:type2est} shows that if $X$ is a UMD space with type $2$ (or, more generally, a Banach space with
martingale type $2$), then for all $p\in [2,\infty)$
the mapping $G\mapsto J(\delta)G$  extends to a bounded operator from
$L^{p}_\mathscr{F}(\R_+\times\O ;\g(H,X))$ to $L^{p}(\R_+\times\O ;X))$
and
the family
\begin{align}\label{eq:calJ}\mathscr{J} := \{J(\delta):\delta>0\}\end{align}
is uniformly bounded.
It what follows, it will be important to know under what additional
conditions this family is $R$-{\em bounded}.

\subsection{$R$-boundedness}

Let $X$ and $Y$ be Banach spaces and let $(r_n)_{n\ge 1}$ be a Radem\-acher
sequence.
A family $\mathscr{T}$ of bounded linear operators from $X$ to $Y$
is called {\em $R$-bounded} if there exists a constant $C\ge 0$ such
that for all finite sequences $(x_n)_{n=1}^N$ in $X$ and
$(T_n)_{n=1}^N$ in ${\mathscr {T}}$ we have
\[ \E \Big\n \sum_{n=1}^N r_n T_n x_n\Big\n^2
\le C^2\E \Big\n \sum_{n=1}^N r_n x_n\Big\n^2.
\]
The least admissible constant $C$ is called the {\em $R$-bound} of
$\mathscr {T}$, notation $R(\mathscr{T})$.  For Hilbert spaces $X$ and $Y$,
 $R$-boundedness is
equivalent to uniform boundedness and $R(\mathscr{T}) = \sup_{t\in
\mathscr{T}} \|T\|$. The notion of
$R$-boundedness has played an important role in recent
progress in the regularity theory of (deterministic) parabolic evolution
equations (see Theorem \ref{thm:maxregdet} below). For more information on
$R$-boundedness and its applications
we refer the reader to \cite{CPSW, DHP,KuWe}.

\medskip

In Theorems \ref{thm:maxregstoch}, \ref{thm:SE}, \ref{thm:SE2}, and \ref{thm:SElocal} it will be important
to have conditions under which the operator family
$\mathscr{J}$ introduced in \eqref{eq:calJ} is not just uniformly
bounded,
but even $R$-bounded, from $L^p_\F(\R_+\times\O ;\g(H,X))$ to $L^p(\R_+\times\O ;X)$.
Whether or not this happens depends on the choice of $p$ and
the geometry of the Banach space $X$. The proof of next proposition
(\cite[Theorem 3.1]{NVW10}) depends critically upon the two-sided estimate provided by
Theorem \ref{thm:UMD}.

\begin{theorem}[Conditions for $R$-boundedness of  $\mathscr{J}$]\label{thm:Rbddness}
In each of the two cases below, $\mathscr{J}$ is $R$-bounded as a family of
operators from $L^p_\F(\R_+\times\O ;\g(H,X))$ to $L^p(\R_+\times\O ;X)$:
\begin{enumerate}[(1)]
\item $p\in [2, \infty)$ and $X$ is isomorphic to
a Hilbert space.
\item $p\in (2, \infty)$ and $X$ is
isomorphic to a closed subspace of $L^q(\OO,\mu)$,
with $q\in (2, \infty)$ and $(\OO,\mu)$ a $\sigma$-finite measure space.
\end{enumerate}
\end{theorem}

The proof of this theorem generalises to $2$-convex UMD Banach lattices $X$ with type $2$
over $(\mathcal{O},\mu)$
whose $2$-concavification $X_{(2)}$ (see \cite[Section 1.d]{LiTz}) is a UMD Banach lattice as well.
Further results about the $R$-boundedness of $\mathscr{J}$ will be contained in a forthcoming
paper \cite{NVW11}.

\section{$H^\infty$-calculi and (stochastic) maximal $L^p$-regularity}

Let $A$ be a sectorial operator, or equivalently, let $-A$ be the generator of a bounded analytic $C_0$-semigroup
$S = (S(t))_{t\ge 0}$ of bounded linear operators on a Banach space $X$.
As is well known (see \cite[Proposition I.1.4.1]{Am}),
the spectrum of $A$ is contained in the closure of a sector
$$\Sigma_{\vartheta} := \{z\in \C\setminus\{0\}:
 |\arg(z)|<\vartheta\}$$ for some $\vartheta\in (0,\frac12\pi)$,
and for all $\sigma\in (\vartheta,\pi)$ one has
\begin{align}\label{eq:parabolic}
\sup_{z\in \C\setminus\Sigma_{\sigma}} \| z(z-A)^{-1}\|<\infty.
\end{align}
In the converse direction, this property characterize negative generators of
bounded analytic $C_0$-semigroups. We refer to \cite{EN, Pa} for more proofs and
further results.

For $\alpha\in (0,1)$ we write
$$
\Dom_A(\alpha,p) =X_{\alpha, p} =  (X, \Dom(A))_{\alpha,p}, \quad
X_{\alpha} = [X, \Dom(A)]_{\alpha}
$$
for the real and complex interpolation scales associated with $A$.
If $A$ has bounded imaginary powers, then (see \cite[Theorem 6.6.9]{Haase:2}, \cite[Theorem 1.15.3]{Tr1})
\begin{align}\label{eq:fracpow}
 X_{\alpha} = \Dom(A^\alpha) \ \hbox{ with equivalent norms}.
\end{align}
The results of Section \ref{sec:4} and Subsection \ref{sec:5a} are of isomorphic nature and
the choice of the norm on $X_\a$
is immaterial. In Subsection \ref{subsec:compH} we shall present a sharp result which is of
isometric nature, for which it
is important to work with the homogeneous norm on $X_\a$ (assuming
bounded invertibility of $A$). We return
to this point in Subsection \ref{subsec:compH}.

We will need the following result (see \cite[Theorem 1.14.5]{Tr1}).

\begin{proposition}\label{prop:initialvalue}
$-A$ be the generator of a bounded analytic $C_0$-semigroup
$S = (S(t))_{t\ge 0}$ of bounded linear operators on a Banach space $X$,
and suppose that $0\in \varrho(A)$. For $x\in X$ the following assertions are
equivalent:
\begin{enumerate}[(1)]
\item The orbit $t\mapsto S(t) x$ belongs to $H^{1,p}(\R_+;X)\cap
L^p(\R_+;\Dom(A))$.
\item The vector $x$ belongs to $\XAp $.
\end{enumerate}
If these equivalent conditions hold, then for all $x\in \XAp $ one has
\[ \max\big\{\|t\mapsto S(t) x\|_{H^{1,p}(\R_+;X)}, \
\|t\mapsto S(t) x\|_{L^p(\R_+;\Dom(A))}\big\}
\eqsim \|x\|_{ \XAp }.\]
\end{proposition}

\subsection{Operators with bounded $H^\infty$-calculus\label{sec:sectorial}}
Let $H^\infty(\Sigma_\sigma)$ denote the
Banach space of all bounded analytic functions
$\varphi:\Sigma_\sigma\to \C$ endowed with the supremum norm.
Let $H_0^\infty(\Sigma_\sigma)$
be its linear subspace consisting of all functions satisfying an estimate
\[
|\varphi(z)|\leq \frac{C|z|^\e}{(1+|z|^2)^\e}
\]
for some $\e>0$.

Now let $-A$ be as above and define,
for $\varphi \in H_0^\infty(\Sigma_\sigma)$ and $\sigma<\sigma'<\pi$,
\[\varphi(A) = \frac{1}{2\pi i}\int_{\partial \Sigma_{\sigma'}} \varphi(z)
(z-A)^{-1}
\, dz. \]
This integral converges absolutely and is independent of $\sigma'$.
We say that $A$ has a {\em bounded $H^\infty(\Sigma_{\sigma})$-calculus}
if there is a constant $C\ge 0$ such that
\begin{equation}\label{eq:Hinfty}
\|\varphi(A)\|\leq C \|\varphi\|_\infty \quad \forall\varphi\in
H^\infty_0(\Sigma_\sigma).
\end{equation}
The least constant $C$ for which this holds will be referred to as the
{\em boundedness constant} of the $H^\infty(\Sigma_{\sigma})$-calculus.
By approximation, the estimate \eqref{eq:Hinfty} can be extended to all
functions $f\in H^\infty(\Sigma_{\sigma})$. The infimum of all $\sigma$ such
that $A$ admits a bounded
$H^\infty(\Sigma_{\sigma})$-calculus is called the {\em angle} of the
calculus.

Any operator $A$ with a bounded $H^\infty$-calculus of angle less than $\frac12\pi$
had bounded imaginary powers. In particular, \eqref{eq:fracpow} applies to such operators.

We proceed with some examples of operators $-A$ for which $A$ has a bounded $H^\infty$-calculus
of angle $<\frac12\pi$; we refer to \cite{DHP, KuWe, Weis-survey}
for further references.
\begin{example}\label{ex:Hinfty} \
\begin{enumerate}
\item Generators of analytic $C_0$-contraction semigroups on Hilbert spaces \cite{McI}.

\item Generators of bounded analytic $C_0$-semigroups admitting Gaussian bounds
\cite{DuRob}.

\item Generators of positive analytic $C_0$-contraction semigroups on a space
$L^q(\mu)$, $1<q<\infty$ \cite{KWcalc}.

\item Second order uniformly elliptic operators \cite{AHS,DDHPV} on
$L^q(\R^d)$ and on $L^q(\OO)$
for bounded $C^{2}$-domains $\OO\subseteq\R^d$ (with Dirichlet or Neumann boundary conditions) \cite{AHS,DDHPV}.

\item The Stokes operator associated with the Navier-Stokes equation
on bounded domains \cite{KKW, NoSa}
(see Section \ref{sec:NS}) and on unbounded
domains \cite{Ku08}.

\item 
Suppose
$-A$ generates a symmetric submarkovian $C_0$-semigroup $S$
on a space $L^2(\mu)$. Then, for all $q\in (1,\infty)$, $A$ admits
a bounded $H^\infty$-calculus of angle $<\frac12\pi$ on $L^q(\mu)$ \cite{KuSt}.

\item The Laplace-Beltrami operator $-A := \Delta_{\rm LB}$ on a complete Riemannian manifold $M$
is given by the symmetric Dirichlet form
$ -\lb \Delta_{\rm LB}f,g\rb = \int_M \nabla f \cdot \nabla g$
and therefore it satisfies the assumptions of example (6)
\cite{Bak87, Str83}.

\item Let $\gamma$ denote the standard Gaussian measure on $\R^n$.
The Ornstein-Uhlen\-beck operator $-A = \Delta_{\rm OU} := \Delta - x\cdot \nabla$
on satisfies the assumptions of example (6). This example admits
various generalisations; see
\cite{ChMGol, Sh92} (for the infinite-dimensional symmetric case) \cite{MPRS} (for the
finite-dimensional non-symmetric case) and \cite{MaaNee} (for the infinite-dimensional
non-symmetric case).

\end{enumerate}
\end{example}

In example (4), under mild assumptions of the coefficients one typically has
$$\Dom(A^\frac12) = H^{1,q}(\R^d) \ \hbox{ and }
 \ H^{1,q}_{\rm Dir/Neum}(\OO)$$ respectively (see, e.g., \cite[Proposition 3.1.7]{Haase:2}
and the references in Sections \ref{sec:parabRd} and \ref{sec:8}).
If, in example (7), the Ricci curvature of $M$ is bounded below, then
$$\Dom((-\Delta_{\rm LB})^\frac12) = H^{1,q}(M),$$ the first order Sobolev space
associated with the derivative $\nabla$ \cite{Bak87}. In example (8),
the classical Meyer inequalities imply that
$$\Dom((-\Delta_{\rm OU})^\frac12) = {\mathbb D}^{1,q}(\R^n,\gamma),$$ the first order Sobolev
space associated with the Malliavin derivative in $L^q(\R^n;\gamma)$ \cite{Nualart}.
Necessary and sufficient conditions for the validity of the analogous identification in
the non-symmetric and infinite-dimensionsional case were obtained in \cite{MaaNee};
special cases were obtained earlier in
\cite{ChMGol, MPRS, Sh92}.

\subsection{Maximal $L^p$-regularity}
Let $-A$ be the generator of a bounded analytic $C_0$-semigroup
$S$ on a Banach space $X$.
For functions $g\in L_{\rm loc}^1(\R_+;X)$
we consider the linear inhomogeneous problem
\begin{equation}\label{eq:detmaxreg}
\left\{
  \begin{aligned}
    u'(t) + A u(t) & = g(t), \qquad  t>0, \\
    u(0) & = 0.
  \end{aligned}
\right.
\end{equation}
The (unique) mild solution to \eqref{eq:detmaxreg} is given by
\[u(t) = S*g (t) := \int_0^t S(t-s) \, g(s) \, d s.\]
Let $p\in (1,\infty)$. For functions $g\in L^p(\R_+;X)$, a
routine estimate shows that for all $\delta\in [0,1)$,
$S*g$ takes values in $\Dom(A^\delta)$ almost everywhere on $\R_+$.
The operator $A$ has {\em maximal $L^p$-regularity} if
for all $g\in L^p(\R_+;X)$ the mild solution $u$
belongs to $\Dom(A)$ almost everywhere on $\R_+$,
and satisfies
\begin{align}\label{eq:MR}
\|A u\|_{L^p(\R_+;X)} \leq C \|g\|_{L^p(\R_+;X)},
\end{align}
where $C$ is a constant independent of $g$.
If $A$ has maximal $L^p$-regularity, then the mild solution $u$
satisfies the identity
$$ u(t) = u_0+ \int_0^t Au(s)\,ds + \int_0^t g(s)\,ds,$$
and the Lebesgue differentiation theorem shows that $u$ is differentiable
almost everywhere on $\R_+$ with derivative $u'(t) = Au(t)+g(t)$.
As a consequence, the inequality \eqref{eq:MR} self-improves to
\begin{align}\label{eq:MR2}
\|u'\|_{L^{p}(\R_+;X)} +
\|A u\|_{L^p(\R_+;X)} \leq C \|g\|_{L^p(\R_+;X)},
\end{align}
with a possibly different constant $C$.

In the definition of maximal $L^p$-regularity we do not insist that $u$
itself be in $L^p(\R_+;X)$.
If, however, $0\in \varrho(A)$, then $Au\in L^p(\R_+;X)$ implies
$u\in L^p(\R_+;X)$, and the estimate \eqref{eq:MR2} is then equivalent to
\begin{align*}
\|u\|_{H^{1,p}(\R_+;X)} + \|u\|_{L^p(\R_+;\Dom(A))} \leq C \|g\|_{L^p(\R_+;X)}.
\end{align*}

The following result was proved in \cite{We} (part (1)) and \cite{KaWe} (part (2));
the final assertion follows by standard trace and interpolation
techniques (see \cite[Theorem III.4.10.2]{Am}).

\begin{theorem}\label{thm:maxregdet}
Let $-A$ be the generator of an analytic $C_0$-semigroup on a UMD space $X$.

\begin{enumerate}[(1)]
 \item The operator $A$ has a maximal $L^p$-regularity  for some (equivalently, all) $p\in (1,\infty)$
if and only if the set
$\{\l (\l+A)^{-1}: \ \l \in i\R\setminus\{0\}\}$ is $R$-bounded in $\calL(X)$.
\item If $A$ has a bounded $H^\infty$-calculus of angle $<\frac12\pi$, then
$A$ has maximal $L^p$-regularity for all $p\in (1,\infty)$.
\end{enumerate}
If $A$ has maximal $L^p$-regularity and $0\in \varrho(A)$, then the mild solution
$u = S*g$ of \eqref{eq:detmaxreg}
belongs to $BUC(\R_+;\XAp )$ and
$$
\|u\|_{BUC(\R_+;\XAp )}\leq C \|g\|_{L^p(\R_+;X)}
$$
with a constant $C$ independent of $g$.
\end{theorem}

\subsection{Stochastic maximal $L^p$-regularity}
In this section we assume that $-A$ generates a bounded analytic
$C_0$-semigroup on a UMD space $X$ with type $2$.
For processes $G\in L^p_{\F}(\R_+\times\O ;\g(H,X))$ we consider the problem
\begin{equation*}
\left\{
  \begin{aligned}
    d U(t) + A U(t) \, dt & = G(t) \, d W_H(t), \qquad t>0, \\
      U(0) & = 0.
  \end{aligned}
\right.
\end{equation*}
The (unique) mild solution of this problem is given by
\[U(t) = \int_0^t S(t-s) \, G(s) \, d W_H(s).\]
Note that this stochastic integral is well defined in view of \eqref{eq:type2est}
and the remark following it. A routine estimate
based on \eqref{eq:type2est} and Young's inequality
shows that for all $\delta\in [0,\frac12)$,
$U$ takes values in $\Dom(A^\delta)$ almost everywhere on $\R_+\times\O$.
The operator $A$ is said to have {\em stochastic maximal $L^p$-regularity} if
for all $G\in L^p_{\F}(\R_+\times\O ;\g(H,X))$, $U$ belongs to $\Dom(A^{\frac12})$ almost everywhere
on $\R_+\times\O $ and
satisfies
\begin{align}\label{eq:stochmaxreg}
\|A^{\frac12} U\|_{L^p(\R_+\times\O ;X)} \leq C
\|G\|_{L^p_{\F}(\R_+\times\O ;\g(H,X))}.
\end{align}
with a constant $C$ independent of $G$.
Under the additional assumption $0\in \varrho(A)$,
$A^{\frac12} U\in {L^p(\R_+\times\O ;X)}$ implies $ U\in
{L^p(\R_+\times\O ;X)}$
and \eqref{eq:stochmaxreg} is equivalent to
\begin{align}\label{eq:stochmaxreg2}
\|U\|_{L^p(\R_+\times\O ;\Dom(A^\frac12))} \leq C
\|G\|_{L^p_{\F}(\R_+\times\O ;\g(H,X))}.
\end{align}

\begin{remark}
It follows from  \cite{NVW10}
that $A$ has stochastic maximal $L^p$-maximal regularity
if and only if \eqref{eq:stochmaxreg} holds for all
deterministic $G\in L^p(\R_+;\gamma(H,X))$.
For later use we note that by Theorem \ref{thm:UMD}, this condition is equivalent to
\begin{align}\label{eq:pathwiseG}
\int_0^\infty \n s\mapsto A^{\frac12} S(t-s) G(s) \|_{\gamma(L^2(0,t;H),X)}^p\,dt
\le C^p\n G\n_{L^p(\R_+;\gamma(H,X))}^p.
\end{align}
\end{remark}

Comparing the notions of deterministic maximal $L^p$-regularity
and stochastic maximal $L^p$-regularity, we note that the latter
increases the regularity only by an exponent $\frac12$.
Another difference is that stochastic maximal $L^p$-regularity does not in general
imply $u\in H^{\frac12,p}(\R_+;L^p(\O;X))$
(see, however, \eqref{eq:stochmaxregmixed} for a related result which does hold true).
In fact (this corresponds to the case $H=X=\R$, $A=0$, and $G$ constant),
already Brownian motions fail to belong to
$H^{\frac12,p}(0,1;L^p(\O))$ for any $p\in [1,\infty]$.
This follows from the continuous inclusion
\[H^{\frac12,p}(0,1;L^p(\O)) \simeq L^p(\O; H^{\frac12,p}(0,1)) \embed L^p(\O; B_{p,p\wedge 2}^\frac12(0,1))\]
and the results in \cite{CKR93, HytVer08}.

Recall the operator family $\mathscr{J}$ which has been introduced
in \eqref{eq:calJ}. By Theorem \ref{thm:Rbddness}, the
$R$-boundedness of $\mathscr{J}$ is satisfied if $X$ is isomorphic to a closed subspace
of an $L^q$-space.

The next theorem has been proved in
\cite[Theorems 1.1, 1.2]{NVW10} for spaces $X=L^q(\mu)$ with $q\ge 2$ and $\mu$ $\sigma$-finite.
Inspection of the proof
shows that it consists of two parts: (i) the proof that $\mathscr{J}$ is $R$-bounded
for such $X = L^q(\mu)$ (\cite[Theorem 3.1]{NVW10},
recalled here as Theorem \ref{thm:Rbddness}) and (ii) the proof that, still for $X=L^q(\mu)$,
the $R$-boundedness of $\mathscr{J}$ implies the result.
Step (ii) extends {\em mutatis mutandis} to arbitrary
UMD Banach space with type $2$, provided one replaces spaces of square functions
such as $L^q(\mu;\mathscr{H})$ and duality
for Hilbert spaces $\mathscr{H}$ by spaces of radonifying operators $\gamma(H,X)$ and trace duality
following the lines of \cite{KaWe}. This leads to the following result:

\begin{theorem}[Conditions for stochastic maximal $L^p$-regularity] \label{thm:maxregstoch}
Let $X$ be a UMD space with type $2$ and let $p\in [2, \infty)$, and suppose the
operator family
$\mathscr{J}$ is $R$-bounded from
$\calL(L^{p}_\mathscr{F}(\R_+\times\O ;\g(H,X))$ to
$L^{p}(\R_+\times\O ;X)$.
If $A$ has a bounded $H^\infty$-calculus on $X$ of angle $<\frac12\pi$,
then $A$ has stochastic maximal $L^p$-regularity.
If, in addition, $0\in \varrho(A)$,
then also \eqref{eq:stochmaxreg2} holds
and
\begin{align}\label{eq:stochmaxregtrace}
\|U\|_{L^p(\O; BUC(\R_+;\Dom_A(\frac12-\frac1p,p))} \leq C
\|G\|_{L^p_{\F}(\R_+\times\O ;\g(H,X))}.
\end{align}
and, for all
$\theta\in [0,\frac12)$,
\begin{align}\label{eq:stochmaxregmixed}
\|U\|_{L^p(\O;H^{\theta,p}(\R_+;\Dom(A^{\frac{1}{2}-\theta})))} \leq C
\, \|G\|_{L^p(\R_+\times\O;\g(H,X))}.
\end{align}
In all these estimates, the constants $C$ are independent of $G$.
\end{theorem}

Note that the case $\theta=0$ of \eqref{eq:stochmaxregmixed}
corresponds to the stochastic maximal $L^p$-regularity estimate \eqref{eq:stochmaxreg}.
The proof of  \eqref{eq:stochmaxregmixed} proceeds by reducing the problem,
via the $H^\infty$-calculus
of $A$, to the $R$-boundedness of a certain family $\mathscr{I}$ of stochastic convolution
operators with scalar-valued kernels. By convexity arguments, the $R$-boundedness of
$\mathscr{I}$ is then deduced from the $R$-boundedness of $\mathscr{J}$.
The estimate
\eqref{eq:stochmaxregtrace} follows from a combination of
\eqref{eq:stochmaxreg}, \eqref{eq:stochmaxregmixed}, and an interpolation
argument (see \cite{Zacher05}).
Note that \eqref{eq:stochmaxregmixed} implies the space-time
H\"older regularity estimate
\begin{align*}
\|U\|_{L^p(\O;C^{\theta-\frac1p}([0,\infty);\Dom(A^{\frac{1}{2}-\theta})))} \leq C
\, \|G\|_{L^p(\R_+\times\O;\g(H,X))}, \qquad \theta\in(\tfrac1p, \tfrac12).
\end{align*}

It has already been observed that the limiting case $\theta= \frac12$ is not allowed in
\eqref{eq:stochmaxregmixed} even when $A=0$
and $G\in\g(H,X)$ is constant.

\section{The main result}\label{sec:4}

On a Banach space $X_0$ we consider the stochastic evolution equation
\begin{equation}\tag{SE}\label{SE}
\left\{\begin{aligned}
dU(t)  +  A U(t)\, dt& = [F(t,U(t)) + f(t)] \,dt \\ & \qquad \qquad + [B(t,U(t)) + b(t)]\,dW_H(t), \qquad t\in
[0,T],\\
 U(0) & = u_0.
\end{aligned}
\right.
\end{equation}

Concerning the space $X_0$, the random operator $A$, the nonlinearities $F$ and
$B$, the external forces $f$ and $b$, and the random initial value $u_0$ we shall assume the following standing hypothesis.

\medskip\noindent
{\bf Hypothesis (H).}

\let\ALTERWERTA\theenumi
\let\ALTERWERTB\labelenumi
\def\theenumi{{\rm (HX)}}
\def\labelenumi{(HX)}
\begin{enumerate}
\item\label{as:XUMD}
$X_0$ is a UMD Banach space with type $2$, and
$X_1$ is a Banach space continuously and densely embedded in $X_0$.
\end{enumerate}
\let\theenumi\ALTERWERTA
\let\labelenumi\ALTERWERTB
\let\ALTERWERTA\theenumi
\let\ALTERWERTB\labelenumi
\def\theenumi{{\rm (HA)}}
\def\labelenumi{(HA)}
\begin{enumerate}
\item\label{as:A}
The function $A:\O\to \calL(X_1, X_0)$ is strongly
$\F_0$-measurable. There exists $w\in\R$ such that each operator $w+A(\omega)$,
viewed as a densely defined operator
on $X_0$ with domain $X_1$,
has a bounded $H^\infty$-calculus of angle
$0< \sigma<\frac12\pi$, with $\sigma$ independent of $\omega$.
There is a constant $C$, independent of $\omega$, such that for all $\varphi\in
H^\infty(\Sigma_{\sigma})$,
\[\|\varphi(w+A(\omega))\|\leq C \|\varphi\|_{H^\infty(\Sigma_{\sigma})}.\]
In what follows, for $\alpha\in (0,1)$ we write
$$
X_{\alpha, p} = (X_0,X_1)_{\alpha,p}, \quad
X_{\alpha} = [X_0,X_1]_{\alpha}
$$
for the real and complex interpolation scales of the couple $(X_0,X_1)$.\end{enumerate}
\let\theenumi\ALTERWERTA
\let\labelenumi\ALTERWERTB

\let\ALTERWERTA\theenumi
\let\ALTERWERTB\labelenumi
\def\theenumi{{\rm (HF)}}
\def\labelenumi{(HF)}
\begin{enumerate}
\item\label{as:LipschitzF}
The function $f:[0,T]\times\O\to X_0$ is adapted and strongly measurable and $f\in L^1(0,T;X_0)$ almost surely.
The function $F:[0,T]\times\O\times X_1\to X_0$ is strongly
measurable and
\begin{enumerate}
\item for all $t\in [0,T]$ and $x\in X_1$ the random variable $\omega\mapsto F(t,\omega,
x)$ is strongly
$\F_t$-measurable;
\item there exist constants
$L_{F}$, $\tilde L_{F}$, $C_{F}$ such that for all $t\in [0,T]$,
$\omega\in \O$, and $x,y\in X_1$,
\begin{equation*}
\phantom{aaaaa}
\|F(t,\omega, x) - F(t,\omega,y)\|_{X_0} \leq L_{F} \|x-y\|_{X_1} +
\tilde L_{F}
\|x-y\|_{X_0}
\end{equation*}
and
\begin{equation*}
\phantom{aa}
\|F(t,\omega, x)\|_{X_0} \leq C_{F}(1+ \|x\|_{X_1}).
\end{equation*}
\end{enumerate}
\end{enumerate}
\let\theenumi\ALTERWERTA
\let\labelenumi\ALTERWERTB

\let\ALTERWERTA\theenumi
\let\ALTERWERTB\labelenumi
\def\theenumi{{\rm (HB)}}
\def\labelenumi{(HB)}
\begin{enumerate}
\item\label{as:LipschitzB}
The function $b:[0,T]\times\O\to \g(H,X_{\frac12})$ is adapted and strongly measurable
and $b\in L^2(0,T;\g(H,X_{\frac12}))$ almost surely.
The function $B:[0,T]\times\O\times X_1\to \g(H,X_{\frac12})$
is strongly measurable and
\begin{enumerate}
\item for all $t\in [0,T]$ and $x\in X_1$ the random variable  $\omega\mapsto
B(t,\omega, x)$ is strongly $\F_t$-measurable;
\item there exist constants $L_{B}$, $\tilde L_{B}$, $C_{B}$ such that for all $t\in
[0,T]$, $\omega\in \O$, and $x,y\in X_1$,
\begin{equation*}
\phantom{aaaaa}
\|B(t,\omega, x) - B(t,\omega,y)\|_{\g(H,X_{\frac12})} \leq L_{B} \|x-y\|_{X_1} +
\tilde L_{B} \|x-y\|_{X_0}
\end{equation*}
and
\begin{equation*}
\phantom{aa}
\|B(t,\omega, x)\|_{\g(H,X_{\frac12})} \leq C_{B}(1+ \|x\|_{X_1}).
\end{equation*}
\end{enumerate}
\end{enumerate}
\let\theenumi\ALTERWERTA
\let\labelenumi\ALTERWERTB

\let\ALTERWERTA\theenumi
\let\ALTERWERTB\labelenumi
\def\theenumi{{\rm (H$u_0$)}}
\def\labelenumi{(H$u_0$)}
\begin{enumerate}
\item \label{as:initial_value}
The initial value $u_0:\O\to X_0$ is strongly $\F_0$-measurable.
\end{enumerate}
\let\theenumi\ALTERWERTA
\let\labelenumi\ALTERWERTB

\begin{remark} \label{rem:ass}
Some comments on these assumptions are in order.
\begin{itemize}

\item[(i)] By \ref{as:A}, the spaces $X_0$ and $X_1$ are isomorphic as Banach
spaces,
an isomorphism being given by $(\lambda-A(\omega))^{-1}$ for any
$\lambda\in\varrho(A(\omega))$.
In particular, since $X_0$ is a UMD space with type $2$,
the same is true for $X_1$. As a consequence, also
the real and complex interpolation spaces $X_{\a,p}$ with $p\in [2, \infty)$ and
$X_\a$ are UMD spaces
with type $2$ (see \cite[Proposition 5.1]{KaMoSm}).

\item[(ii)] If \ref{as:A} holds for some $w\in\R$, then it holds for
any $w'>w.$ Furthermore, we may write $$-A + F = -(A +w') + (F + w'),$$
and note that a function $F$ satisfies \ref{as:LipschitzF} if and only if
$F+w'$ satisfies \ref{as:LipschitzF}.
Thus, in what follows we may replace $A$ and $F$ by $A+w'$ and $F+w'$ and thereby
assume, without any loss of generality,
that the operators $A(\omega)$ are {\em invertible}, uniformly
in $\omega$.

\item[(iii)] The operators
$-A(\omega)$ generate analytic $C_0$-semigroups $S(\omega)$
on $X_0$, given through the $H^\infty$-calculus by $$S(t,\omega) =
e^{-tA(\omega)},\quad t\ge 0.$$
For each $t\ge 0$ and $x\in X_0$,
$\omega\mapsto S(t,\omega)x$ is strongly $\F_0$-measurable.
Assuming, as in (ii), that the operators $A(\omega)$ are uniformly invertible,
the semigroups $S(\cdot,\omega)$ are uniformly exponentially stable, uniformly in
$\omega$.

\item[(iv)] By \eqref{eq:pathwiseG}, Theorem \ref{thm:maxregstoch} extends to the present situation
of a random operator $A$ satisfying \ref{as:A}.

\item[(v)]
The Lipschitz conditions in \ref{as:LipschitzF} and \ref{as:LipschitzB}
are fulfilled if and only if there exist $\a_F,\a_B\in [0,1)$
and constants $L_F', \tilde L_F', L_B', \tilde L_B'$ such that
\[
\phantom{aaaaa}
\|F(t,\omega, x) - F(t,\omega,y)\|_{X_0} \leq L_{F}' \|x-y\|_{X_1} +
\tilde L_{F}' \|x-y\|_{X_{\a_F}}
\]
and
\[
\phantom{aaaaa}
\|B(t,\omega, x) - B(t,\omega,y)\|_{\g(H,X_{\frac12})} \leq L_{B}' \|x-y\|_{X_1} +
\tilde L_{B}' \|x-y\|_{X_{\a_B}}.
\]
Moreover, for any $\e>0$ the constants
$\tilde L_F'$ and $\tilde L_B'$ can be chosen in such a
way that $|L_F' - L_F|<\e$ and $|L_B' - L_B|<\e$.
The `if' part is obvious from $\|x-y\|_{X_0} \lesssim_\alpha \|x-y\|_{X_\a}$
(in this case we may take $L_F' = L_F $ and $L_B' = L_B$),
and the `only
if' part follows by a standard application of Young's inequality.
Indeed, for any $\delta>0$ we have
\[\phantom{xxx}\|x-y\|_{X_\a}\leq C \|x-y\|^{1-\alpha}_{X_0} \|x-y\|^{\alpha}_{X_1}\leq \frac{C}{(1-\alpha)\delta} \|x-y\|_{X_0} + \frac{C\delta}{\alpha} \|x-y\|_{X_1}.\]
Choosing $\delta>0$ small enough this gives the required result.
In certain applications (see Sections \ref{sec:parabRd}, \ref{sec:8} and \ref{sec:NS} below) this reformulation
of the conditions \ref{as:LipschitzF} and \ref{as:LipschitzB} is more convenient.
\end{itemize}
\end{remark}

\begin{definition}\label{def:strongsol}
Let {\rm (H)}   be satisfied.
A process $U: [0,T]\times\Omega \to X_0$ is called a
{\em strong solution} of \eqref{SE} if it is strongly measurable and adapted,
and
\begin{enumerate}[(i)]
\item almost surely, $U\in L^2(0,T;X_1)$;

\item for all $t\in [0,T]$, almost surely the following identity holds in $X_0$:
\begin{align*}
\phantom{aa}
U(t) + \int_0^t A U(s) \, ds = u_0 & + \int_0^t F(s,U(s)) + f(s) \, ds \\ & + \int_0^t
B(s,U(s)) +b(s) \, d W_H(s).
\end{align*}
\end{enumerate}
\end{definition}

To see that the integrals in this definition are well defined, we
note that, by \ref{as:A}, the process $A U$ is strongly measurable and satisfies
\[\|A U\|_{L^1(0,T;X_0)} \leq \|A\|_{\calL(X_1,X_0)} \|U\|_{L^1(0,T;X_1)}\]
almost surely. Similarly, by \ref{as:LipschitzF} and
\ref{as:LipschitzB}, $F(\cdot, U)$ and $ f$ belong to $ L^1(0,T;X_0)$ and $B(\cdot, U)$ and $b$
belong to $L^2(0,T;\g(H,X_{\frac12}))$ almost surely. Therefore, the Bochner integral is
well defined in $X_0$, and the stochastic integral is well defined in
$X_{\frac12}$ (and hence in $X_0$) by \ref{as:XUMD},
the fact the space $X_{\frac12}$ is a UMD space with type $2$, and \eqref{eq:type2est}.

By Definition \ref{def:strongsol}, a strong
solution always has a version with continuous paths in $X_0$ such that,
almost surely, the identity in (ii) holds for all $t\in [0,T]$. Indeed, define
$\tilde{U}:[0,T]\times\O\to X_0$ by
\begin{align*}
\tilde{U}(t) := -\int_0^t A U(s) \, ds +  u_0 & + \int_0^t F(s,U(s)) +f(s)\, ds
\\ & + \int_0^t B(s,U(s)) + b(s) \, d W_H(s),
\end{align*}
where we take continuous versions of the integrals on the right-hand side. From
the definitions of $U$ and $\tilde U$ one obtains, for all $t\in [0,T]$,
that $U(t) =
\tilde{U}(t)$ almost surely in $X_0$.
Therefore, almost surely, for all $t\in [0,T]$ one has
\begin{align*}
\tilde{U}(t) +\int_0^t A \tilde{U}(s) \, ds =  u_0 & + \int_0^t F(s,\tilde{U}(s)) + f(s)
\, ds \\ &  + \int_0^t B(s,\tilde{U}(s)) + b(s) \, d W_H(s).
\end{align*}
From now on we choose this version whenever this is convenient.
We will actually prove much stronger regularity properties in Theorem
\ref{thm:SE} below.

\begin{definition}\label{def:mildsol}
Let {\rm (H)} be satisfied.
A process $U: [0,T]\times\Omega \to X_0$ is called a
{\em mild solution} of \eqref{SE} if
it is strongly measurable and adapted, and
\begin{enumerate}[(i)]
\item almost surely, $U\in L^2(0,T;X_1)$;

\item  for all $t\in [0,T]$, almost surely
the following identity holds in $X_0$:
\begin{align*}
U(t) = S(t) u_0 & + \int_0^t S(t-s)[ F(s,U(s)) + f(s)]\,ds
\\ & + \int_0^t S(t-s)[B(s,U(s))+ b(s)]\,dW_H(s).
\end{align*}
\end{enumerate}
\end{definition}

The convolutions with $F(\cdot,U(\cdot))$ and $f$ are well defined as an $X_0$-valued
process by
\ref{as:LipschitzF}. The stochastic convolutions with $B(\cdot,U(\cdot))$ and $b$ are
well defined as
an $X_{\frac12}$-valued process
(and hence as an $X_0$-valued process) by \ref{as:LipschitzB}, the fact that
$X_{\frac12}$ is a UMD space with type $2$,
and \eqref{eq:type2est}.
Henceforth we shall use the notations
\begin{align*}
S*g(t) &:= \int_0^t S(t-s)g(s)\,ds, \\
S\diamond G(t) &:= \int_0^t S(t-s)G(s)\,dW_H(s),
\end{align*}
whenever the integrals are well defined.

\begin{proposition}\label{prop:strongmild}
Let {\rm (H)}   be satisfied.
A process $U:[0,T]\times\Omega\to X_0$ is a strong solution
of \eqref{SE} if and only if it is a mild solution of \eqref{SE}.
\end{proposition}
Results of this type for time-dependent operators $A$ are well known.
Since in our case $A$ also depends on $\O$, the usual
proof has to be adjusted. For the reader's convenience we provide the details.
\begin{proof}
For notational convenience we write $\overline{F}(t,x) = F(t,x) + f(t)$ and $\overline{B}(t,x) = B(t,x) + b(t)$.

First assume that $U$ is a mild solution. As in
\cite[Proposition 6.4 (i)]{DPZ}, the (stochastic) Fubini theorem can be used to
show
that for all $t\in [0,T]$, almost surely we have
\begin{equation*}
\begin{aligned}
U(t) + \int_0^t AU(s)\, ds  =  u_0
 + \int_0^t \overline{F}(s,U(s)) \, ds + \int_0^t \overline{B}(s,U(s)) \, d W_{H}(s).
\end{aligned}
\end{equation*}

Next assume that $U$ is a strong solution of \eqref{SE}.
By the scalar-valued It\^o formula,
\begin{equation*}
\begin{aligned}
\lb  U(t), \varphi(t)\rb - \lb u_0, \varphi(0)\rb
& =\int_0^t \lb AU(s),\varphi(s)\rb + \lb U(s), \varphi'(s)\rb\, ds
\\ & \qquad +
\int_0^t \lb \overline{F}(s,U(s)), \varphi(s)\rb \, ds
\\ & \qquad +
\int_0^t \overline{B}(s,U(s))^* \varphi(s) \, d W_{H}(s),
\end{aligned}
\end{equation*}
for functions $\varphi\in C^1([0,t];E^*)$ of the form $\varphi= g\otimes x^*$.
By linearity and density this extends to all $\varphi\in
C^1([0,t];E^*)$. By linearity and approximation this extends to all $\varphi\in
L^0(\O;C^1([0,t];E^*))$ which are $\F_0$-measurable. Indeed, recall that for a
$\limn \int_0^\cdot \psi(t) - \psi_n(t) \, d W(t)= 0$ in
$L^0(\O;C([0,T]))$ whenever $\limn \psi_n =\psi$ in $L^0(\O;L^2(0,T;H))$ (see
\cite[Proposition 17.6]{Kal}).

With the choice
$\varphi(t)=S^*(t-s) \lambda (\lambda+A^*)^{-1} x^*$
we obtain, for all $x^*\in E^*$ and
$\lambda>w$ (with $w$ as in \ref{as:A}),
\begin{align*}
\lb  \lambda(\lambda+A)^{-1}  U(t), &x^*\rb - \lb \lambda(\lambda+A)^{-1} S(t)
u_0,  x^* \rb
\\ & =\Big<\lambda(\lambda+A)^{-1} \int_0^t \lb S(t-s) \overline{F}(s,U(s)) \, ds,  x^*
\Big>
\\ & \qquad + \Big<\lambda(\lambda+A)^{-1}  \int_0^t S(t-s)\overline{B}(s,U(s))  \, d W_{H}(s),
x^*\Big>,
\end{align*}
where we used the strong $\F_0$-measurability of $A$.
Now the result follows from the fact that for all $x\in X$,
$\lambda(\lambda+A)^{-1} x \to x$ as $\lambda\to \infty$.
\end{proof}

Let us fix an exponent $p\in [2,\infty)$ for the moment and assume, as in Remark \ref{rem:ass}(ii),
that the operators $A(\omega)$ are uniformly invertible.
By Theorem \ref{thm:maxregdet} and \ref{as:XUMD} and \ref{as:A}, the linear operator
$$g\mapsto S*g$$ is bounded from
$L^p_{\F}(\O;L^p(\R_+;X_0))$ into $L^p_{\F}(\O;L^p(\R_+;X_1))$.
Furthermore, if
the operator family $\mathscr{J}$ introduced in Subsection \ref{sec:J}
is $R$-bounded in
\[\calL(L^{p}_\mathscr{F}(\R_+\times\O ;\g(H,X_{0})),
L^{p}(\R_+\times\O ;X_{0})),\]
then it is also $R$-bounded in
\[\calL(L^{p}_\mathscr{F}(\R_+\times\O ;\g(H,X_{\frac12})),
L^{p}(\R_+\times\O ;X_{\frac12}))\]
and therefore by Theorem \ref{thm:maxregstoch} (applied to the space
$X_{\frac12}$) and \ref{as:XUMD} and \ref{as:A}, the reiteration identity $X_1 =
(X_{\frac12})_{\frac12}$ (apply $A^{\frac12}$ to both sides
and use that $\Dom(A^\frac12 ) = X_\frac12$ by \ref{as:A})
the mapping
$$G\mapsto S\diamond G$$ is bounded from
$L^p_{\F}(\O;L^p(\R_+;\gamma(H,X_{\frac12})))$
into $L^p_{\F}(\O;L^p(\R_+;X_1))$.
We shall denote by $$K^*_p \ \hbox{ and } \  K^\diamond_p$$
the norms of these operators. We emphasise that the numerical value of these
constants depends on the choice of the parameter $w'$ used for rescaling $A$
(cf. remark \ref{rem:ass}).

In what follows we fix an arbitrary time horizon $T>0$; constants
appearing in the inequalities below are allowed to depend on it.
Recall that by Theorem \ref{thm:Rbddness}, the
$R$-boundedness of the operator family
$\mathscr{J}$ is satisfied if $X_0$ is isomorphic to a closed subspace
of an $L^q$-space.

\begin{theorem}\label{thm:SE}
Let {\rm (H)} be satisfied, let $p\in [2,\infty)$,
let $f\in L^p_{\F}(\O;L^p(0,T;X_0))$ and $b\in L^p_{\F}(\O;\allowbreak L^p(0,T;\allowbreak
\g(H,X_\frac12)))$,
and suppose that the operator family
$\mathscr{J}$ is $R$-bounded from
$\calL(L^{p}_\mathscr{F}(\R_+\times\O ;\g(H,X_{0}))$ to
$L^{p}(\R_+\times\O ;X_0)$. If the Lipschitz constants $L_F$ and $L_B$
satisfy
$$K^*_p L_F  + K^\diamond_p L_B <1,$$
then the following assertions hold:
\begin{enumerate}[(i)]
\item
If $u_0\in L_{\F_0}^0(\O;\Xp)$,  then the problem \eqref{SE}
has a unique strong solution $U$ in
$$L^0_{\F}(\O;\allowbreak L^p(0,T;X_1)) \cap L^0_{\F}(\O;\allowbreak C([0,T];\Xp )).$$
\item If $u_0\in L_{\F_0}^p(\Omega; \Xp )$, then the strong solution $U$ given by part {\rm (i)}
belongs to $$L^p_{\F}((0,T)\times \O;X_1)\cap
L_{\F}^p(\O;C([0,T];\Xp ))$$  and satisfies
\begin{align*}
\|U\|_{L^p((0,T)\times\O;X_1)} &\leq
C(1+\|u_0\|_{L^p(\O;\Xp )}),\\
\phantom{aaaa} \|U\|_{L^p(\O;C([0,T];\Xp ))} & \leq
C(1+\|u_0\|_{L^p(\O;\Xp )}),
\end{align*}
with constants $C$ independent of $u_0$.

\item For all $u_0, v_0\in L_{\F_0}^p(\O;\Xp )$,
the corresponding strong solutions $U,V$ satisfy
\begin{align*}
\|U-V\|_{L^p((0,T)\times\O;X_1)} &\leq C
\|u_0-v_0\|_{L^p(\O;\Xp )},\\
\phantom{aaaa} \|U-V\|_{L^p(\O;C([0,T];\Xp ))} & \leq C\|u_0-
v_0\|_{L^p(\O;\Xp )},
\end{align*}
with constants $C$ independent of $u_0$ and $v_0$.
\end{enumerate}
\end{theorem}

\begin{remark}
The condition $u_0\in L_{\F_0}^0(\Omega;\Xp )$ is
satisfied if \ref{as:initial_value} holds and
$u_0$ takes values in $\Xp $ almost surely.
Indeed, by \ref{as:initial_value} we know that $u_0$
is strongly $\F_0$-measurable as an $X$-valued random variable.
Now the  strong $\F_0$-measurability of $u_0$ as an $\Xp $-valued
random variable easily follows from the strong measurability of $\xi:\Omega\to
L^p(0,1,\frac{dt}{t};X)$, given by
$$ \xi(\omega):= [t\mapsto  AS(t) u_0(\omega)],$$
and the definition of $\Xp $.
\end{remark}

\begin{proof}[Proof of Theorem \ref{thm:SE}]
Without loss of generality we can reduce to the case where $w=0$
(see Remark \ref{rem:ass} (ii)).
By assumption we have
$K^*_p L_F  + K^\diamond_p L_B = 1-\theta$ for some $\theta\in (0,1]$.
Without loss of generality we may assume that $L_F + L_B >0$ and $\theta\in
(0,1)$.

By Proposition \ref{prop:strongmild}
it suffices to prove existence and uniqueness of a mild solution.

\smallskip
{\em Step 1:} Local existence of mild solutions for initial values $u_0\in
L_{\F_0}^p(\O;\Xp )$.
We fix a number $\kappa \in (0,T]$, to be chosen in a moment, and
introduce, for $\theta\in [0,1]$, the Banach spaces
\begin{align*}
Z_{\theta,\kappa} & = L^p_{\F}(\O;L^p(0,\kappa;X_\theta)),\\
Z^\gamma_{\theta,\kappa} & = L^p_{\F}(\O;L^p(0,\kappa;\g(H,X_{\theta}))).
\end{align*}
On $Z_{1,\kappa}$ we define an equivalent norm $\nnn\cdot\nnn$ by
\[\nnn\phi\nnn = \|\phi\|_{Z_{1,\kappa}} + M  \|\phi\|_{Z_{0,\kappa}}\]
with $M = (K^*_p \tilde L_{F} + K^\diamond_p \tilde L_{B})/(K^*_p
L_F + K^\diamond_p L_B)$.

In order to simplify notations we shall omit the subscript $\kappa$ in what
follows.
Let $L:Z_1\to Z_1$ be the mapping given by
\[L(\phi)(t) = S(t) u_0 + S*[F(\cdot, \phi)+f](t) + S\diamond [B(\cdot, \phi)+b](t).\]
We emphasise that $L$ depends on the initial value $u_0$.

First we check that $L$ does indeed map $Z_1$ into itself.
By \ref{as:initial_value} and Proposition \ref{prop:initialvalue},
$t \mapsto S(t)u_0$ defines an element of $Z_1$.

By restriction to the interval $[0,\kappa]$, the operators
$g\mapsto S*g$ and $G\mapsto S\diamond G$ are bounded as mappings from
$L^p_{\F}(\O;L^p(0,\kappa;X_0))$ and
$L^p_{\F}(\O;L^p(0,\kappa;\gamma(H,X_{\frac12})))$
into $L^p_{\F}(\O;L^p(0,\kappa;X_1))$, with norms bounded by
$K^*_p$ and $K^\diamond_p$ respectively.
Therefore $L$ is well defined as a mapping from $Z_1$ into itself,
and for all $\phi_1,\phi_2\in Z_1$ we may estimate
\begin{align*}
\|L(\phi_1)-L(\phi_2)\|_{Z_1} &\leq \|S*(F(\cdot,
\phi_1)-F(\cdot, \phi_2))\|_{Z_1} + \|S\diamond (B(\cdot, \phi_1) -
B(\cdot, \phi_2))\|_{Z_1}
\\ & \leq K^*_p \|F(\cdot, \phi_1)-F(\cdot,
\phi_2)\|_{Z_0} + K^\diamond_p \|B(\cdot, \phi_1) - B(\cdot,
\phi_2)\|_{Z_{\frac12}^\g}
\\ & \leq K^*_p L_F \|\phi_1-\phi_2\|_{Z_1} + K^*_p \tilde L_{F}
\|\phi_1-\phi_2\|_{Z_0} \\ & \qquad   + K^\diamond_p L_B \|\phi_1 -
\phi_2\|_{Z_1} + K^\diamond_p \tilde L_{B} \|\phi_1 - \phi_2\|_{Z_{0}}
\\ & =  (1-\theta) \nnn \phi_1-\phi_2\nnn,
\end{align*}
recalling that $K^*_p L_F  + K^\diamond_p L_B = 1-\theta$.
Moreover, we have the elementary estimate
\begin{align*}
\|L(\phi_1)-L(\phi_2)\|_{Z_0} & \leq c(\kappa) \big[C L_F
\|\phi_1-\phi_2\|_{Z_1} + C \tilde L_{F} \|\phi_1-\phi_2\|_{Z_0}  \\ &
\qquad\qquad +C' L_B \|\phi_1 - \phi_2\|_{Z_1} + C' \tilde L_{B} \|\phi_1 -
\phi_2\|_{Z_0}\big]
\\ & \leq \tilde{c}(\kappa) \nnn \phi_1 - \phi_2\nnn,
\end{align*}
where
$\kappa\mapsto c(\kappa) $ and $\kappa\mapsto \tilde c(\kappa)$ are continuous
functions on $[0,T]$ not depending on $u_0$ and satisfying
$\lim_{\kappa\downarrow 0} c(\kappa) = \lim_{\kappa\downarrow 0}\tilde c(\kappa) =
0$.

Collecting the above estimates, we see that
\[\nnn L(\phi_1)-L(\phi_2) \nnn \leq (1-\theta + M\tilde{c}(\kappa)) \nnn \phi_1
-
\phi_2\nnn.\]
So far, the number $\kappa>0$ was arbitrary. Now we set
\[\kappa : = \inf\{t\in (0,T]: \, M\tilde{c}(t) \geq \tfrac12\theta\}.\]
where we take $\kappa = T$ if the infimum is taken over the empty set. Note that
$\kappa$ only depends on $\theta$, the Lipschitz constants of $F$ and $B$, the
constants $K^*_p$ and $K^\diamond_p$ and the type $2$ constant of $X_{\frac12}$.
Then $(1-\theta + M\tilde{c}(\kappa))\leq 1-\tfrac12\theta$,
and it follows that $L$
has a unique fixed point in $Z_1$. This gives a
process $U\in Z_1$ such that for almost all $(t,\omega)\in
[0,\kappa]\times\O$, the following identity holds in $X_1$:
\begin{equation}\label{eq:mildfix}
U(t) = S(t) u_0 + S*F(\cdot,U)(t) + S*f(t) + S\diamond B(\cdot,U)(t) + S \diamond b(t).
\end{equation}
By Theorems \ref{thm:maxregdet} and \ref{thm:maxregstoch} (applied with
$X=X_{\frac12}$), and keeping in mind Remarks \ref{rem:ass}(i) and (iv),
 $U$ has a version with trajectories in
$L^p(\O;C([0,\kappa],  \Xp ))$.
For this version, almost surely the identity \eqref{eq:mildfix} holds
in $X_0$ for all $t\in [0,\kappa]$.

\smallskip

{\em Step 2:} Local existence of mild solutions for initial values
$u_0\in
L_{\F_0}^0(\O;\Xp )$.

For $n\geq 1$, let
$$\Gamma_n := \big\{\|u_0\|_{\Xp }\leq n\big\}.
$$
From Step 1 we obtain processes
$U_n$ belonging to $ Z_1\cap L^p(\O;C([0,\kappa], \Xp ))$ such that
\eqref{eq:mildfix} holds with the pair $(u_0, U)$ replaced by
$(u_{0,n},U_n)$ (with $u_{0,n} = \one_{\Gamma_n} u_0$). We claim that for all $m\leq n$, $U_n(\cdot,\omega)
= U_m(\cdot,\omega)$ in $ \Xp $ almost surely on
$\Gamma_m\times [0,\tau_m]$.
Indeed, by Step 1 and the fact that $\Gamma_m\in\F_0$,
\[\begin{aligned}
\nnn\one_{\Gamma_m} (U_m -U_n)\nnn &=  \nnn\one_{\Gamma_m} (L(U_m) -L(U_n))\nnn
\\ &  = \nnn\one_{\Gamma_m}(L(\one_{\Gamma_m} U_m) -L(\one_{\Gamma_m} U_n))\nnn
\\ & \leq \nnn L(\one_{\Gamma_m} U_m) -L(\one_{\Gamma_m} U_n)\nnn
\\ & \leq (1-\tfrac12\theta)\, \nnn\one_{\Gamma_m}(U_m-U_n)\nnn
\end{aligned}\]
and since $\theta\in (0,1)$ it follows that for almost all $(t,\omega)\in
[0,\kappa]\times\Gamma_n$, $U_m(t,\omega) = U_n(t,\omega)$ in $X_1$. By
\eqref{eq:mildfix} for $U_m$ and $U_n$ it follows that for almost
all $\omega\in \Gamma_m$, $U_n(\cdot, \omega) = U_m(\cdot, \omega)$
in $ \Xp $, and the claim follows.
Therefore, we can define
$U:[0,\kappa]\times\O\to X_0$ by $U = U_n$ on $\Gamma_n$. Now it is easy to
check that
\[U\in L^0_{\F}(\O;L^p(0,\kappa;X_1))\cap L^0(\O;C([0,\kappa],
\Xp )).\]
and that for all $t\in [0,\kappa]$, \eqref{eq:mildfix} holds almost surely in
$X_0$.

\smallskip

{\em Step 3:} Local uniqueness of mild solutions for initial values $u_0\in
L_{\F_0}^0(\O;\Xp )$.

Let $U,V\in L^0_{\F}(\O;L^p(0,\kappa;X_1))$ be such that
\eqref{eq:mildfix} holds. For $W\in \{U,V\}$ let $\tau^W_n$ be the
stopping time defined by
\[\tau^W_n = \inf\{t\in [0,\kappa]: \|\one_{[0,t]}W\|_{L^p(0,\kappa;X_1)}\geq
n\}\]
(and $\tau_n^W = \kappa$ if this set is empty)
and let $\tau_n = \tau_n^U\wedge \tau_n^V$. Let $U_n =
\one_{[0,\tau_n]}U$ and $V_n = \one_{[0,\tau_n]}V$. Clearly, for all
$n\geq 1$, we have $U_n, V_n\in Z_1$. Using the extension of \cite[Lemma
A.1]{BMS} to the type $2$ setting one can check that for
all $t\in [0,\kappa]$, almost surely, one has
\begin{align*}
W_n = \one_{[0,\tau_n]} S(\cdot) u_0 & +  \one_{[0,\tau_n]}
(S*(\one_{[0,\tau_n]} (F(\cdot, W_n)+ f)))
\\ & + \one_{[0,\tau_n]} (S\diamond (\one_{[0,\tau_n]} (B(\cdot, W_n)+b)))
\end{align*}
in $X_0$, where $W_n \in\{ U_n,V_n\}$. As in Step 1 it follows that
\[\begin{aligned}
\nnn & U_n - V_n\nnn \\ &  \leq \nnn S*(\one_{[0,\tau_n]} (F(\cdot,
U_n) - F(\cdot,V_n) ))\nnn + \nnn S\diamond (\one_{[0,\tau_n]} (B(\cdot, U_n) -
B(\cdot,V_n) ))\nnn
\\ & \leq (1-\tfrac12\theta) \nnn U_n - V_n \nnn.
\end{aligned}\]
Since $\theta\in (0,1)$, we obtain that $U_n = V_n$ in $Z_1$. Letting $n$
tend to infinity, we may conclude that $U=V$ in
$L^0_{\F}(\O;L^p(0,\kappa;X_1))$.

\smallskip

{\em Step 4:} Global existence of mild solutions.

In Steps $1$ and $2$ we have shown that there exists a unique mild solution
$U_1$ in $L^0_{\F}(\O;L^p(0,\kappa;X_1))$ with trajectories in
$C([0,\kappa],  \Xp )$. Let, for $0\le a<b\le T$,
\[Y(a,b) :=  L^0_{\F}(\O;L^p(a,b;X_1)) \cap L^0(\O;C([a,b], \Xp ))\]
We construct a mild solution on $[\kappa, 2\kappa]$.
Using the path continuity in $ \Xp $ we can take $u_{\kappa}
= U_1(\kappa)$ in $L^0(\O; \Xp )$ as initial value and repeat
Steps $1$ and $2$ to obtain a unique mild solution $U_2\in Y(\kappa, 2\kappa)$
on $[\kappa, 2\kappa]$ with initial data $u_{\kappa}$. One easily
checks that letting $U = U_1$ on $[0,\kappa]$ and $U=U_2$ on $[\kappa, 2\kappa]$
defines a mild solution on $[0,2\kappa]$. Iterating this finitely many times we
obtain a mild solution $U\in Y(0,T)$.

\smallskip
{\em Step 5:} Global uniqueness of mild solutions.

To see that $U$ is the unique mild solution in $Y(0,T)$, let $V$ be
another mild solution in $Y(0,T)$. Recall from Step 1 that we can find versions of $U$ and
$V$ which also have paths in $C([0,T];\Xp )$. It suffices to prove
the uniqueness for these versions.
Note that by the uniqueness on $[0,\kappa]$ we have
$U|_{[0,\kappa]} = V|_{[0,\kappa]}$. By the almost sure pathwise continuity
of
$U$ and $V$ with values in the space $\Xp $ we see that
almost surely $U(\kappa) = V(\kappa)$ in $\Xp $. One easily
checks that both $U|_{[\kappa,2\kappa]}$ and $V|_{[\kappa,2\kappa]}$ are mild
solutions in $Y(\kappa, 2\kappa)$ on the interval $[\kappa, 2\kappa]$. By
uniqueness on
$[\kappa, 2\kappa]$ from Step 3a, we obtain that $U|_{[\kappa,2\kappa]} =
V|_{[\kappa,2\kappa]}$ in $Y(\kappa, 2\kappa)$. Proceeding in finitely many
steps we obtain $U =  V$ in $Y(0,T)$.

\smallskip

{\em Step 6:} The proof of part (ii).

On $[0,\kappa]$  it follows from Step 1 that
\begin{align*}
\nnn U \nnn = \nnn L(U) \nnn & \le \nnn L(U) - L(0)\nnn  + \nnn L(0)\nnn
\\ & \leq (1-\tfrac12\theta) \nnn U\nnn + C(1+\|u_0\|_{L^p(\O;\Xp )}).
\end{align*}
Since $\theta\in (0,1)$ we obtain
\begin{equation}\label{eq:UZ1est}
\nnn U \nnn \leq \frac{2C}{\theta}  (1+\|u_0\|_{L^p(\O;\Xp )}).
\end{equation}

Next, observe that by Proposition \ref{prop:initialvalue}, Theorems
\ref{thm:maxregdet} and \ref{thm:maxregstoch}, Remark \ref{rem:ass} (iv),
and \ref{as:LipschitzF} and \ref{as:LipschitzB} one has
\begin{align*}
\ & \|U\|_{L^p(\O;C([0,\kappa];\Xp ))} \\ &
\qquad =
\|L(U)\|_{L^p(\O;C([0,\kappa];\Xp ))}
\\ & \qquad \leq C \|u_0\|_{L^p(\Omega;\Xp )}
+ K^*_p \|F(\cdot, U)+f\|_{Z_0} + K^\diamond_p\|B(\cdot, U)+b
\|_{Z_{\frac12}^\g}
\\ & \qquad \leq C \|u_0\|_{L^p(\Omega;\Xp )}
+ K^*_p  C_{F,f} (1+\|U\|_{Z_1}) + K^\diamond_p C_{B,b} (1+\|U\|_{Z_1}).
\end{align*}
From \eqref{eq:UZ1est} and the norm equivalence of $\nnn\cdot\nnn$ on $Z_1$
we obtain
\begin{equation}\label{eq:UZ1est2}
\|U\|_{L^p(\O;C([0,\kappa];\Xp ))}\leq \tilde{C}
(1+\|u_0\|_{L^p(\Omega; \Xp )})
\end{equation}
for some constant $\tilde{C}$. This proves the required estimates on $[0,
\kappa]$. In particular, it follows from \eqref{eq:UZ1est2} that
\begin{equation}\label{eq:Ukappa}
\|U(\kappa)\|_{L^p(\Omega;\Xp )}\leq \tilde{C}
(1+\|u_0\|_{L^p(\Omega;\Xp )}).
\end{equation}
Using $U(\kappa)$ as an initial values the same argument now gives the following
estimates for $U$ on $[\kappa, 2\kappa]$:
\begin{align*}
\|U\|_{L^p_{\F}(\O;L^p(\kappa, 2\kappa;X_1))} &\leq \frac{2C}{\theta}
(1+\|U(\kappa)\|_{L^p(\Omega; \Xp )})
\\ \|U\|_{L^p(\O;C([\kappa, 2\kappa];\Xp ))}&\leq \tilde{C}
(1+\|U(\kappa)\|_{L^p(\Omega; \Xp )}).
\end{align*}
Combining this with \eqref{eq:Ukappa} and iterating this finitely many times
gives (2).

\smallskip
{\em Step 7:} The proof of part (iii).

First note that by Step 1,
\[\begin{aligned}
\|U-V\|_{Z_1} & = \|L(U)-L(V)- S u_0 + S v_0\|_{Z_1} \\
&  \leq (1-\tfrac12\theta)\|U-V\|_{Z_1} + C\|u_0-v_0\|_{L^p(\O;\Xp )},
\end{aligned}\]
where $L=L_{u_0}$ is the operator from Step 1 with initial condition $u_0$.

Since $\theta\in (0,1)$ this implies
\[\|U-V\|_{Z_1} \leq \frac{2C}{\theta}\|u_0-v_0\|_{L^p(\O;\Xp )}.\]
In the same way as for \eqref{eq:UZ1est2} one can prove that
\[\|U - V\|_{L^p(\O;C([0,\kappa];\Xp ))}\leq
\tilde{C}\|u_0-v_0\|_{L^p(\O;\Xp )}.\]
Now one iterates the argument as in Steps 4 and 5.
\end{proof}

Theorem \ref{thm:SE} can be seen as an extension of \cite{Brz2} to the borderline
case. A maximal $L^p$-regularity result using real
interpolation spaces instead of fractional domain spaces has been obtained in \cite{Brz1}.

\begin{remark}
We believe that by using Lenglart's inequality (see \cite{Lenglart}),
it may be shown that in Theorem \ref{thm:SE} one obtains solutions in
$L^{p_1}(\O;L^{p_2}(0,T;X_1))$
and $L^{p_1}(\O;C([0,T];X_{1-\frac1p_2,p_2}))$
for any $p_2>p_1>0$ and $p_2\ge 2$.
Since we do not have any applications of this, we shall not pursue this any further.
\end{remark}

\begin{remark}
Applying \eqref{eq:stochmaxregmixed} to the space $X_{\frac12}$ one can prove in the
same way that
\[U\in L^0(\O;H^{\theta,p}(0,T;X_{1-\theta})) \ \text{for all $\theta\in
[0, \tfrac12)$}.\]
In particular,
\[U\in L^0(\O;C^{\theta-\frac1p}([0,T];X_{1-\theta})) \ \text{for all
$\theta\in [\tfrac1p, \tfrac12)$}.\]
Moreover, the following estimates hold:
\begin{align*}
\|U\|_{L^p(\O;H^{\theta,p}(0,T;X_{1-\theta}))} &\leq
C(1+\|u_0\|_{L^p(\O;\Xp )}),
\\ \|U-V\|_{L^p(\O;H^{\theta,p}(0,T;X_{1-\theta}))} &\leq C
\|u_0-v_0\|_{L^p(\O;\Xp )},
\end{align*}
where $U$ and $V$ are the solutions with initial values $u_0$ and $v_0$
respectively.
\end{remark}

\section{Extensions of the main result}

\subsection{The time-dependent case}\label{sec:5a}

In the same setting as before we now consider \eqref{SE} with an adapted
operator family $\{A(t,\omega):\,t\in [0,T], \, \omega\in\O\}$ in $\calL(X_1, X_0)$:
\begin{equation}\tag{SE$'$}\label{SE'}
\left\{\begin{aligned}
dU(t)  +  A(t) U(t)\, dt & =  [F(t,U(t)) + f(t)]\,dt \\ & \qquad  \qquad  + [B(t,U(t)) +b(t)]\,dW_H(t), \qquad t\in
[0,T],\\
 U(0) & = u_0.
\end{aligned}
\right.
\end{equation}
Below we shall extend the definition of a strong solution
(see Definition \ref{def:strongsol})
to the time-dependent problem \eqref{SE'} for adapted random operators
$A:[0,T]\times\O\to \calL(X_1,X_0)$.
There is no direct extension of the definition of a mild solution to this
setting, the reason being that serious problems with adaptedness arise
 (see \cite{LeNu98} for details). Below we shall prove the existence and uniqueness
of strong solutions for \eqref{SE'} by means of maximal regularity techniques.

Throughout this section we replace Hypothesis \ref{as:A}
by the following hypothesis \ref{as:A'} and
we say that {\em Hypothesis {\rm (H)}$'$ holds} if \ref{as:XUMD}, \ref{as:A'}, \ref{as:LipschitzF}, \ref{as:LipschitzB} and \ref{as:initial_value} hold,
with

\let\ALTERWERTA\theenumi
\let\ALTERWERTB\labelenumi
\def\theenumi{{\rm (HA)$'$}}
\def\labelenumi{(HA)$'$}
\begin{enumerate}
\item\label{as:A'}
The function $A:[0,T]\times\O\to \calL(X_1, X_0)$ is strongly measurable and
adapted.
Each operator $A(t,\omega)$,
viewed as a densely defined operator
on $X_0$ with domain $X_1$, is invertible and
has a bounded $H^\infty$-calculus of angle
$0< \sigma<\frac12\pi$, with $\sigma$ independent of $t$ and $\omega$.
There is a constant $C$, independent of $t$ and $\omega$, such that for all $\varphi\in
H^\infty(\Sigma_{\sigma})$,
\[\|\varphi(A(t,\omega))\|\leq C \|\varphi\|_{H^\infty(\Sigma_{\sigma})}.\]
The function $A:[0,T]\times\O\to \calL(X_1,X_0)$ is {\em piecewise
relatively continuous, uniformly in $\omega$}, i.e., there exists finitely many points
$0=t_0<t_1<\ldots<t_{N}=T$ such that for all $\e>0$ there exists a $\delta>0$
and $\eta>0$ such that for all $\omega\in \O$, for all $1\leq n\leq N$, for all
$t,s\in [t_{n-1},t_n]$ and for all $x\in X_1$, we have
\[ \ \qquad |t-s|<\delta \ \implies \ \|A(t,\omega)x - A(s,\omega)x\|_{X_0}<\e\|x\|_{X_1}  + \eta \|x\|_{X_0}.\]
\end{enumerate}
\let\theenumi\ALTERWERTA
\let\labelenumi\ALTERWERTB
The first part of Hypothesis \ref{as:A'} implies that the operators $-A(t,\omega)$ generate
bounded analytic $C_0$-semigroups on $X_0$ for which the estimate \eqref{eq:parabolic}
holds uniformly in $t$ and $\omega$.

Relatively continuous operators $A$ have been introduced in \cite{ACFP} to study
maximal $L^p$-regularity for deterministic problems. We consider a piecewise variant
here, which seems to be new even in a deterministic setting. It seems that
the results in \cite{ACFP} extend to this more general setting without difficulty.

\begin{definition}
Let {\rm (H)$'$} be satisfied.
A process $U: [0,T]\times\Omega \to X_0$ is called a
{\em strong solution} of \eqref{SE'} if it is strongly measurable and adapted,
and
\begin{enumerate}[(i)]
\item almost surely, $U\in L^2(0,T;X_1)$;
\item for all $t\in [0,T]$, almost surely the following identity holds in $X_0$:
\begin{equation}
\label{strongsol'}\begin{aligned}
U(t) + \int_0^t A(s) U(s) \, ds  = u_0 & +  \int_0^t F(s,U(s)) + f(s) \, ds \\ &  + \int_0^t
B(s,U(s)) + b(s)\, d W_H(s).
\end{aligned}
\end{equation}
\end{enumerate}
\end{definition}

As before, under {\rm (H)$'$} all integrals are well defined, and again $U$ has a pathwise continuous version
for which, almost surely, the identity in (ii) holds for all $t\in [0,T]$.

\begin{theorem}\label{thm:SE2}
Let {\rm (H)$'$} be satisfied, let $p\in [2,\infty)$, and suppose that the operator family
$\mathscr{J}$ is $R$-bounded from
$L^{p}_\mathscr{F}(\R_+\times\O ;\g(H,X_{0}))$ to
$L^{p}(\R_+\times\O ;X_0)$.
If the Lipschitz constants $L_F$ and $L_B$
satisfy $$K^*_p L_F  + K^\diamond_p L_B <1, $$
then the assertions of Theorem \ref{thm:SE} (i),
(ii) and (iii) remain true for the problem \eqref{SE'}.
\end{theorem}

\begin{proof}
As in the proof of Theorem \ref{thm:SE} we may assume that
$K^*_p L_F  + K^\diamond_p L_B = 1-\theta$ with $\theta\in (0,1)$.

Choose $\delta>0$ and $\eta>0$ such that for all $1\leq n\leq N$ and for all
$t,s\in [t_{n-1},t_n]$, for all $x\in X_1$,
\[\|A(t) x - A(s)x\|_{X_0}\leq \tfrac12{\theta} \|x\|_{X_1} + \eta \|x\|_{X_0}
\text{ if } |t-s|<\delta.\]
Fix $0=s_0<s_1< \ldots< s_M=T$ such that $\{t_{n}:0\leq n\leq N\}$ is a subset
of $\{s_{n}:0\leq n\leq N\}$ and $|s_{m}-s_{m-1}|<\delta$ for $m=1, \ldots, M$.

We first solve the problem on $[0, s_1]$. Let $F_{A,0}:[0,s_1]\times\O\times
X_1\to X_0$ be defined by $F_{A,0}(t,x) = F(t,x) - A(t)x + A(0)x$. Then $F_A$
satisfies \ref{as:LipschitzF} with $F$ replaced by $F_{A,0}$. Moreover,
$L_{F_{A,0}}\leq L_F + \frac12\theta$ and $\tilde L_{F_{A,0}}\leq
\tilde L_{F} + C \eta$, and therefore, the condition of Theorem \ref{thm:SE}
holds for the equation with $F$ replaced by $F_{A,0}$ and $A$ replaced by $A(0)$
with constant $K^*_p L_{F_{A,0}}  + K^\diamond_p L_B =1-\frac12\theta$. Therefore,
Theorem \ref{thm:SE} implies the existence of a unique strong solution $U\in
L^0_{\F}(\O; L^p(0,s_1;X_1))$, i.e.
almost surely, for all $t\in [0,s_1]$ the following identity holds in $X_0$:
\begin{align*}
U(t) + \int_0^t A(0) U(s) \, ds = u_0& + \int_0^t F_{A,0}(s,U(s)) + f(s) \, ds \\ & + \int_0^t
B(s,U(s)) + b(s) \, d W_H(s)
\end{align*}
and therefore also \eqref{strongsol'} holds on $[0,s_1]$ almost surely.
Moreover, the assertions of Theorem \ref{thm:SE} (i), (ii) and (iii) hold on $[0,s_1]$.

Now we proceed inductively. Suppose we know that the assertions
of Theorem \ref{thm:SE} (i), (ii) and (iii) hold for the problem \eqref{SE'}
on the interval $[0, s_{m}]$ with $m\leq M$. If $m=M$, there is
nothing left to prove. If $m<M$, we shall prove next existence and uniqueness on
the interval $[s_{m},
s_{m+1}]$.

Consider the problem
\begin{equation}\label{eq:VFAn}
\left\{\begin{aligned}
dV(t)  +  A(s_{m}) V(t)\, dt& = [F_{A,{m}}(t,V(t)) + f(t)] \,dt \\ & \qquad \qquad
 + [B(t,V(t))+ b(t)]\,dW_H(t),\quad t\in [s_{m}, s_{m+1}],\\
 V(s_{m-1}) & = U(s_{m})
\end{aligned}
\right.
\end{equation}
with $F_{A,m} = F(t,x) - A(t) + A(s_{m})$. As before, Theorem \ref{thm:SE}
(more precisely, the version of it with initial time $0$ replaced by $s_m$) can be
applied to obtain a unique strong solution $V\in L^0_{\F}(\O;
L^p(s_{m},s_{m+1};X_1))$ and assertions (i), (ii) and (iii) of Theorem \ref{thm:SE} hold
for the solution $V$ of \eqref{eq:VFAn}. Now we extend $U$ to $[0,s_{m+1}]$ by
setting $U(t) := V(t)$ for $t\in [s_m,s_{m+1}]$. Then $U$ is in $L^0_{\F}(\O;
L^p(0,s_{m+1};X_1))$ and, using the induction hypothesis, one sees that it is a
strong solution on $[0,s_{m+1}]$. It is also the unique strong solution on
$[0,s_{m+1}]$. Indeed, let $W\in L^0_{\F}(\O; L^p(0,s_{m+1};X_1))$ be another
strong solution on $[0,s_{m+1}]$. By the induction hypothesis we
have $W = U$ in $L^0_{\F}(L^p(0,s_{m};X_1))$. In particular, the definition of a
strong solution implies that $W(s_{m}) = U(s_m)$ almost surely. Now one can see
that $W$ is strong solution of \eqref{eq:VFAn} on $[s_{m},s_{m+1}]$. Since the
solution of \eqref{eq:VFAn} is unique, it follows that also $W = V$ in
$L^0_{\F}(\O;L^p(s_{m},s_{m+1};X_1))$. Therefore, the definition of $U$ shows
that $U = W$ in $L^0_{\F}(L^p(0,s_{m+1};X_1))$. The other results in (i), (ii)
and (iii) for $U$ on $[0,s_{m+1}]$ follow from the corresponding results for $V$
as well. This completes the induction step and the proof.
\end{proof}

\subsection{The locally Lipschitz case}\label{subsec:local}

In this section we shall prove an extension of Theorem
\ref{thm:SE} to the case where the functions $F$ and $B$ satisfy a local Lipschitz
condition with respect to the $\Xp$-norm, where $p\in [2, \infty)$ is fixed.
We replace the Hypotheses \ref{as:LipschitzF} and replace \ref{as:LipschitzB}
by the hypotheses \ref{as:LipschitzFlocal} and
\ref{as:LipschitzBlocal}.

\medskip\noindent
{\bf Hypothesis (H)$_{\text{loc}}^p$}

\let\ALTERWERTA\theenumi
\let\ALTERWERTB\labelenumi
\def\theenumi{{\rm (HF)$_{\rm{loc}}^p$}}
\def\labelenumi{(HF)$_{\rm{loc}}^p$}
\begin{enumerate}
\item\label{as:LipschitzFlocal}
The function $f:[0,T]\times\O\to X_0$ is adapted and strongly measurable and $f\in L^1(0,T;X_0)$ almost surely.
The function $F:[0,T]\times\O\times X_1\to X_0$ is given by
$F = F^{(1)} + F^{(2)}$, where $F^{(1)}:[0,T]\times\O\times X_1\to X_0$
and $F^{(2)}:[0,T]\times\O\times \Xp\to X_0$ are strongly measurable.
The function $F^{(1)}$ is $\F$-adapted and Lipschitz continuous, i.e., it
satisfies \ref{as:LipschitzF}:
\begin{enumerate}[(a)]
\item for all $t\in [0,T]$ and $x\in X_1$ the random variable $\omega\mapsto F^{(1)}(t,\omega,
x)$ is strongly $\F_t$-measurable;
\item there exist constants
$L_{F^{(1)}}$, $\tilde L_{F^{(1)}}$, $C_{F^{(1)}}$ such that for all $t\in [0,T]$,
$\omega\in \O$, and $x,y\in X_1$,
\begin{equation*}
\phantom{aa}
\|F^{(1)}(t,\omega, x) - F^{(1)}(t,\omega,y)\|_{X_0} \leq L_{F^{(1)}} \|x-y\|_{X_1} +
\tilde L_{F^{(1)}} \|x-y\|_{X_0}
\end{equation*}
and
\begin{equation*}
\phantom{aaa}
\|F^{(1)}(t,\omega, x)\|_{X_0} \leq C_{F^{(1)}}(1+ \|x\|_{X_1}).
\end{equation*}
\end{enumerate}
The function $F^{(2)}$ is $\F$-adapted and locally Lipschitz continuous, i.e.,
\begin{enumerate}[(a)]
\setcounter{enumii}{2}
\item for all $t\in [0,T]$ and $x\in \Xp$ the random variable $\omega\mapsto F^{(2)}(t,\omega,
x)$ is strongly $\F_t$-measurable;
\item for all $R>0$ a constant $L_{F^{(2)},R}$ such that for all $t\in [0,T]$, $\omega\in \O$, and $x,y\in X_1$ satisfying $\n x\n_{\Xp}, \n y\n_{\Xp} \le R$,
\begin{equation*}
\phantom{aaaaa}
\|F^{(2)}(t,\omega, x) - F^{(2)}(t,\omega,y)\|_{X_0} \leq L_{F^{(2)},R} \|x-y\|_{\Xp}
\end{equation*}
and there exists a constant $C_{F^{(2)}}$ such that for all $t\in [0,T]$, $\omega\in \O$,
\begin{equation*}
\phantom{aa}
\|F^{(2)}(t,\omega, 0)\|_{X_0} \leq C_{F^{(2)}}.
\end{equation*}
\end{enumerate}
\end{enumerate}
\let\theenumi\ALTERWERTA
\let\labelenumi\ALTERWERTB

\let\ALTERWERTA\theenumi
\let\ALTERWERTB\labelenumi
\def\theenumi{{\rm (HB)$_{\rm{loc}}^p$}}
\def\labelenumi{(HB)$_{\rm{loc}}^p$}
\begin{enumerate}
\item\label{as:LipschitzBlocal}
The function $b:[0,T]\times\O\to \g(H,X_{\frac12})$ is adapted and strongly measurable and $b\in L^2(0,T;\g(H,X_{\frac12}))$ almost surely.
The function $B:[0,T]\times\O\times X_1\to \g(H,X_{\frac{1}{2}})$ is given by
$B = B^{(1)} + B^{(2)}$, where $B^{(1)}:[0,T]\times\O\times X_1\to
\g(H,X_{\frac{1}{2}})$ and $B^{(2)}:[0,T]\times\O\times \Xp\to
\g(H,X_{\frac{1}{2}})$ are strongly measurable. The function $B^{(1)}$
is $\F$-adapted and Lipschitz continuous, i.e., it satisfies
\ref{as:LipschitzB}:
\begin{enumerate}[(a)]
\item for all $t\in [0,T]$ and $x\in X_1$ the random variable $\omega\mapsto B^{(1)}(t,\omega,
x)$ is strongly $\F_t$-measurable;
\item there exist constants
$L_{B^{(1)}}$, $\tilde L_{B^{(1)}}$, $C_{B^{(1)}}$ such that for all $t\in [0,T]$,
$\omega\in \O$, and $x,y\in X_1$,
\begin{equation*}
\phantom{aa}
\|B^{(1)}(t,\omega, x) - B^{(1)}(t,\omega,y)\|_{\g(H,X_{\frac{1}{2}})} \leq L_{B^{(1)}} \|x-y\|_{X_1} +
\tilde L_{B^{(1)}} \|x-y\|_{X_0}
\end{equation*}
and
\begin{equation*}
\phantom{a}
\|B^{(1)}(t,\omega, x)\|_{\g(H,X_{\frac{1}{2}})} \leq C_{B^{(1)}}(1+ \|x\|_{X_1}).
\end{equation*}
\end{enumerate}
The function $B^{(2)}$ is $\F$-adapted and locally Lipschitz continuous, i.e.,
\begin{enumerate}[(a)]
\setcounter{enumii}{2}
\item for all $t\in [0,T]$ and $x\in \Xp$ the random variable $\omega\mapsto B^{(2)}(t,\omega,x)$ is strongly $\F_t$-measurable;
\item for all $R>0$ a constant $L_{B^{(2)},R}$ such that for all $t\in [0,T]$, $\omega\in \O$, and $x,y\in X_1$ satisfying $\n x\n_{\Xp}, \n y\n_{\Xp} \le R$,
\begin{equation*}
\phantom{aaaaa}
\|B^{(2)}(t,\omega, x) - B^{(2)}(t,\omega,y)\|_{\g(H,X_{\frac{1}{2}})} \leq L_{B^{(2)},R} \|x-y\|_{\Xp}
\end{equation*}
and there exists a constant $C_{B^{(2)}}$ such that for all $t\in [0,T]$, $\omega\in \O$,
\begin{equation*}
\phantom{a}
\|B^{(2)}(t,\omega, 0)\|_{\g(H,X_{\frac{1}{2}})} \leq C_{B^{(2)}}.
\end{equation*}
\end{enumerate}
\end{enumerate}
\let\theenumi\ALTERWERTA
\let\labelenumi\ALTERWERTB

Before we explain the definition of a local mild solution, we need to discuss some preliminaries on stopped stochastic convolutions.
Let $G:[0,T]\times\O\to \g(H,X_0)$ be an adapted process which satisfies $G\in L^2(0,T;\g(H,X_0))$ almost surely. Let $\tau$ be a stopping time with values in $[0,T]$.
Define the $X_0$-valued processes $I(G)$ by
\begin{align*}
I(G)(t) = \int_0^t S(t-s) G(s)\, d W_H(s).
\end{align*}
As explained in \cite{BMS} it is tempting to write
\[I(G)(t\wedge \tau) = \int_0^{t\wedge \tau}  S(t\wedge \tau-s) G(s)\, d W_H(s).\]
This is meaningless, however, since the integrand in the right-hand expression
is not adapted, and therefore the stochastic integral is not well defined.
To remedy this problem, following \cite{BMS} we consider the
process $I_{\tau}(G)$ defined by
\[I_{\tau}(G)(t) = \int_0^t \one_{[0,\tau]}(s) S(t-s) G(s)\, d W_H(s)
 = S \diamond (\one_{[0,\tau]}G).\]

The following lemma can be proved as in \cite[Lemma A.1]{BMS}.
\begin{lemma}\label{lem:BMS}
Assume \ref{as:XUMD}.
Let $G:[0,T]\times\O\to \g(H,X_0)$ be an adapted process which satisfies
$G\in L^2(0,T;\g(H,X_0))$ almost surely. Let $\tau$ be a stopping time with
values in $[0,T]$.
If the processes $I(G)$ and $I_{\tau}(G)$ have an $X_0$-valued continuous version,
then almost surely,
\[S(t-t\wedge \tau) I(G)(t) = I_{\tau}(G)(t), \ \ t\in [0,T].\]
In particular, almost surely,
\[I(G)(t\wedge \tau) = I_{\tau}(G)(t\wedge \tau), \ \ t\in [0,T].\]
\end{lemma}
Note that if $G$ is only defined up to a stopping time $\tau'$
with $\tau\leq \tau'$ and $ \one_{[0,\tau]}G $ is in $L^2(0,T;\g(H,X_0))$,
the above definition of $I_{\tau}(G)$ is still meaningful.
This is what we will use below.

\begin{remark}
If \ref{as:A} holds
and $G$ belongs to $L^p(0,T;\g(H,X_0))$
almost surely for some $p>2$, then Theorem \ref{thm:maxregstoch}
(combined with Remark \ref{rem:ass} (iv)) shows that $I(G)$ and
$I_{\tau}(G)$ are both pathwise continuous as $X_{\frac12-\frac1p,p}$-valued
processes, hence also as $X_0$-valued processes.
For $p=2$, pathwise continuity of $I(G)$ and $I_{\tau}(G)$ follows from
\cite[Theorem 1.1]{VW11}.
\end{remark}

\begin{definition}\label{def:mildsolloc}
Let $p\in [2, \infty)$ and let {\rm(H)$_{\rm{loc}}^p$} be satisfied.
Let $\tau$ be a stopping time with values in $[0,T]$. A process $U: [0,\tau)\times\Omega \to \Xp$ is called a
{\em local mild solution} of \eqref{SE} if
$U$ is adapted and for each $\omega\in \O$, $t\mapsto U(t,\omega)$
is continuous in $\Xp$ on the interval $[0,\tau(\omega))$ and, for all $n\ge 1$,
\begin{enumerate}[(i)]
\item almost surely, $\one_{[0,\tau_n]}U\in L^2(0,T;X_1)$;

\item almost surely, for all $t\in [0,T]$, the following identity holds in $X_0$:
\begin{align*}
U(t\wedge \tau_n) = S(t\wedge \tau_n) u_0 & + \int_0^{t\wedge \tau_n} S(t\wedge \tau_n-s)[F(s,U(s)) + f(s)]\,ds
\\ & + S\diamond (\one_{[0,\tau_n]}(B(\cdot,U) + b)))(t\wedge \tau_n),
\end{align*}
where
\[\tau_n = \inf\{t\in [0,\tau): \|U(t)\|_{\Xp}\geq n\}.\]
\end{enumerate}
\end{definition}
Note that
\[S\diamond (\one_{[0,\tau_n]}(B(\cdot,U) + b)))(t\wedge \tau_n) = I_{\tau_n}(B(\cdot,U) + b)(t\wedge \tau_n).\]
The motivation for this expression has been explained in Lemma \ref{lem:BMS}.

A process $U: [0,\tau)\times\Omega \to \Xp$ is called a
{\em maximal local mild solution on $[0,T]$}
if it is a local mild solution and for every stopping time
$\tau'$ with values in $[0,T]$ and every local mild solution
$V:[0,\tau')\times\O\to \Xp$ one has $\tau = \tau'$ almost surely
and $U= V$ in $C([0,\tau);\Xp)$ almost surely.
A process $U: [0,T)\times\Omega \to \Xp$ is called a
{\em global mild solution} if $U$ is a local mild solution (with
$\tau = T$) and $U\in L^2(0,T;X_1)$ almost surely. For such $U$ one
easily checks that part (ii) of Definition \ref{def:mildsol} holds.

In a similar way one can define local and global strong solutions.
It is obvious from the proof of Proposition \ref{prop:strongmild} that
the notions of global strong solution and global mild solution are equivalent.
Below we shall only consider local and global mild solutions.

The following theorem can be proved by following the lines of \cite{Brz1, Sei} (see also  \cite{BMS} and
\cite[Theorem 8.1]{NVW3}).

\begin{theorem}\label{thm:SElocal}
Let {\rm (H)$_{\text{loc}}^p$} be satisfied for $p\in [2,\infty)$, and suppose that the operator family
$\mathscr{J}$ is $R$-bounded from
$L^{p}_\mathscr{F}(\R_+\times\O ;\g(H,X_{0}))$ to $L^{p}(\R_+\times\O ;X_0)$. If the Lipschitz constants $L_F$ and $L_B$
satisfy
$$K^*_p L_F  + K^\diamond_p L_B <1,$$
then the following assertion holds:
\begin{enumerate}[(i)]
\item
If $u_0\in L_{\F_0}^0(\O;\Xp)$, $f\in L_{\F_0}^0(\O;L^p(0,T;X_0))$, and
$b\in L_{\F_0}^0(\O;L^p(0,T;\allowbreak \g(H,X_\frac12)))$, then the problem \eqref{SE}
has a unique maximal local mild solution $U$ in
$$L^0_{\F}(\O; \allowbreak L^p(0,\tau;X_1)) \cap L^0_{\F}(\O; C([0,\tau); \allowbreak \Xp )).$$
\item
If, in addition to the assumptions in (i), $F^{(2)}$ and $B^{(2)}$ also satisfy the linear growth conditions
\begin{align*}
\phantom{aaa}
\|F^{(2)}(t,\omega, x)\|_{X_0} & \leq C_{F^{(2)}}(1+ \|x\|_{\Xp}), \\
\|B^{(2)}(t,\omega, x)\|_{\g(H,X_{\frac12})} & \leq C_{B^{(2)}}(1+ \|x\|_{\Xp}),
\end{align*}
for some constants $C_{F^{(2)}}$ and $C_{B^{(2)}}$ independent of $t\in [0,T]$, $\omega\in \O$, and  $x\in \Xp$,
then the solution $U$ in (i) is a global mild solution which belongs to
$$L^0_{\F}(\O;\allowbreak L^p(0,T;X_1)) \cap L^0_{\F}(\O;\allowbreak C([0,T];\Xp )).$$
\item If, in addition to the assumptions of (i) and (ii), we have $u_0\in L_{\F_0}^p(\Omega; \Xp )$,
$f\in L^p_{\F}(\O;L^p(0,T;X_0))$, and $b\in L^p_{\F}(\O;L^p(0,T;\allowbreak
\g(H,X_\frac12)))$, then the global solution $U$ in (ii)
belongs to $L^p_{\F}((0,T)\times \O;X_1)\cap
L_{\F}^p(\O;C([0,T];\Xp ))$  and satisfies
\begin{align*}
\|U\|_{L^p((0,T)\times\O;X_1)} &\leq
C(1+\|u_0\|_{L^p(\O;\Xp )}),\\
\phantom{aaaa} \|U\|_{L^p(\O;C([0,T];\Xp ))} & \leq
C(1+\|u_0\|_{L^p(\O;\Xp )}),
\end{align*}
with constants $C$ independent of $u_0$.
\end{enumerate}
\end{theorem}

\subsection{The Hilbert space case\label{subsec:compH}}

For Hilbert spaces $X_0$, several of the constants in the estimates in Theorems
\ref{thm:SE} and \ref{thm:SE2} become explicit and we can give more
precise conditions on the smallness of $L_F$ and $L_B$.
Below, we show that if $A$ is self-adjoint and positive,
then $K^*_2 \leq 1$ and $K^{\diamond}_2\leq \frac1{\sqrt 2}$
(these constants have been defined in the text preceding Theorem \ref{thm:SE}).
Moreover, these estimates are optimal in the sense that the condition \eqref{eq:sharpL}
below cannot be improved (see \cite[Section 4.0]{Rozov} for the stochastic part; see also
\cite{BrzVer11} for more information on the smallness condition for
$K^*_p$ and $K^{\diamond}_p$ for $p\not=2$).
As a consequence one obtains the
following result, which is well known to experts (see \cite{DPZ,Rozov} for
related results and \cite{DS04} for applications to a class of SDPEs).

\begin{corollary}
Let $X_0$ and $X_1$ be Hilbert spaces, and let $A:[0,T]\times\O\to
\calL(X_1,X_0)$ be strongly measurable, adapted, self-adjoint,
and piecewise relatively continuous uniformly on $\Omega$. Moreover,
assume that there is a
constant $\delta >0$ such that
\[\|e^{s A(t,\omega)}\|\leq e^{-\delta s},  \ \ t\in [0,T], \ \omega\in\O.\]
Assume \ref{as:LipschitzF}, \ref{as:LipschitzB} and
\ref{as:initial_value}. The assertions of Theorem \ref{thm:SE2}
hold whenever
\begin{align}\label{eq:sharpL}
L_F  + \frac{L_B}{\sqrt{2}}  <1.
\end{align}
\end{corollary}

A similar consequence of Theorem \ref{thm:SElocal} can be formulated in the Hilbert space setting.
\begin{proof}
The result follows at once from Theorem \ref{thm:SE2}
once we show that  $K^*_2 \leq 1$ and $K^{\diamond}_2\leq \frac1{\sqrt 2}$.
Here is it important to endow $X_{\frac12}$ with the norm
\begin{align}\label{eq:honogeneous} \n x\n_{\frac12} := \n A^\frac12 x\n
\end{align}
(cf. the discussion below \eqref{eq:fracpow}).
By the invertibility of $A$ and the equivalence of norms \eqref{eq:fracpow},
\eqref{eq:honogeneous} indeed defines an equivalent norm on $X_{\frac12}$. If what follows, we understand
$K^*_2$ and $K^{\diamond}_2$ as the operator norms as defined in Section \ref{sec:4},
with $X_\frac12$ normed by \eqref{eq:honogeneous}.

We first show that $K^*_2 \leq 1$. Using the spectral theorem one can see that
for all $s\in \R$, one has
\begin{align}\label{eq:spectral}\|A(is+A)^{-1}\|\leq 1\end{align}
As direct proof is obtained as follows. For $x\in X_0$ with $\|x\|\leq1$ and $s\in \R$ one has
\begin{align*}
\|A (is+A)^{-1}x\|^2 & = \lb A^2 (-is+A)^{-1} (is+A)^{-1}x, x \rb \\ & = \lb
A^2 (s^2+A^2)^{-1}x, x \rb = \lb A^2 (t+A^2)^{-1}x, x \rb=:f(t),
\end{align*}
where $t=s^2$ Then $f(0) = 1$ and, for $t>0$,
\[f'(t) = - \lb A^2 (t+A^2)^{-2}x, x \rb = - \|A (t+A^2)^{-1}x\|^2\leq
0,\]
and therefore $f(t)\le 1$ as claimed.

By \eqref{eq:spectral} and Plancherel's theorem, for any $g\in L^2(\R_+;X_0)$
one has that
\begin{align*}
\|A S* g\|_{L^2(\R_+;X_0)}^2 & = \int_\R \|A(is+A)^{-1}\hat{g}(s)\|_{X_0}^2 \,
ds \leq \int_\R \|\hat{g}(s)\|_{X_0}^2 \, ds  = \|g\|_{L^2(\R_+;X_0)}^2,
\end{align*}
and hence $K^*_2\leq 1$.

Next we show that $K^{\diamond}_2\leq \frac1{\sqrt 2}$ (cf. \cite[Section 6.3.2]{DPZ}).
By standard arguments involving the essentially separable-valuedness of strongly
measurable mappings (cf. \cite{NeeCMA}) there is no loss of generality in assuming that
that $H$ is separable. Let $(h_n)_{n\ge 1}$ be an orthonormal basis of $H$.
Let $\calL_2(H,X_\frac12)$ denote the space of Hilbert-Schmidt operators (which is
canonically isometric to $\g(H,X_\frac12)$).
By the It\^o isometry, for all $G\in L^2(\R_+\times\O ;\calL_2(H,X_\frac12))$
we have
\begin{align*}
\|A^{\frac12} S\diamond G\|_{L^2(\R_+\times\O ;X_\frac12)}^2 & =  \int_0^\infty \int_0^t
\sum_{n\geq 1} \E  \|AS(t-s) G(s) h_n \|^2 \, d s \, dt
\\ & \leq  \int_0^\infty \int_0^\infty \sum_{n\geq 1} \E  \|AS(t) G(s)
h_n \|^2 \, dt \, d s
\\ & = \sum_{n\geq 1} \E   \int_0^\infty \int_0^\infty [A^2 S(2t) G(s) h_n, G(s)
h_n] \, d t \, ds
\\ & = \sum_{n\geq 1} \E   \int_0^\infty \tfrac12 [Ag(s) h_n, G(s) h_n] \, ds
\\ & = \tfrac12 \|G\|_{L^2(\R_+\times\Omega;\calL_2(H,X_\frac12))}^2.
\end{align*}
 It follows that $K^{\diamond}_2\leq \frac1{\sqrt 2}$.
\end{proof}

\section{Parabolic SPDEs of order $2m$ on $\R^d$}\label{sec:parabRd}

In this section we shall apply our abstract results to
the following system of $N$ coupled stochastic partial differential
equations on $[0,T]\times \R^d$:
\begin{equation}\label{eq:eq2m}
\left\{
\begin{aligned}
d u(t,x) +  \mathcal{A}(t,x, D) u(t,x)\, dt &= [f(t,x,u) + f^0(t,x)]\,dt \\ & \qquad +
 \sum_{i\geq 1}[ b_i(t,x,u) + b_i^0(t,x) ]\, d w_i(t),
\\ u(0,x) & = u_0(x).
\end{aligned}
\right.
\end{equation}
Here
\[\mathcal{A}(t,\omega,x, D) = \sum_{|\alpha|\leq 2m} a_{\alpha}(t,\omega,x)
D^{\alpha},\]
with $D = -i(\partial_1, \ldots, \partial_d)$.
The precise assumptions on the coefficients
$$a_{\alpha}:[0,T]\times\O\times\R^d\to \C^{N}\times\C^N$$
and the functions
\begin{align*}
f:[0,T]\times\O\times\R^d\times
H^{2m,q}(\R^d;\C^N) & \to L^q(\R^d;\C^N) \\
f^0:[0,T] & \to L^q(\R^d;\C^N)\\
b_i:[0,T]\times\O\times\R^d\times H^{2m,q}(\R^d;\C^N) & \to H^{m,q}(\R^d;\C^N)\\
b^0:[0,T] & \to H^{m,q}(\R^d;\C^N)
\end{align*}
will be stated in the next two subsections. Essentially, we shall assume
that the conditions of \cite{DuSi} (where the non-random case was discussed) hold
pointwise on $\Omega$ with uniform bounds.

\subsection{Hypotheses on the coefficients $a_\alpha$}

Let $\mathcal{A}_\pi$ be the principal part of $\mathcal{A}$,
\[\mathcal{A}_{\pi}(t,\omega,x, D) = \sum_{|\alpha|= 2m}
a_{\alpha}(t,\omega,x) D^\alpha.\]

\let\ALTERWERTA\theenumi
\let\ALTERWERTB\labelenumi
\def\theenumi{{\rm (Ha)}}
\def\labelenumi{(Ha)}
\begin{enumerate}
\item\label{as:Acoef}
The coefficients $a_{\alpha}:[0,T]\times\O\times\R^d\to \C^{N}\times\C^N$
are $\mathcal{P}\times
\mathcal{B}_{\R^d}$-measurable, where $\mathcal{P}$ denotes the progressive
$\sigma$-algebra of $[0,T]\times\Omega$ and $\mathcal{B}_{\R^d}$ the Borel $\sigma$-algebra of ${\R^d}$.
Furthermore, \smallskip

\begin{enumerate}[(i)]
\item
$a_{\alpha} \in L^\infty(\O;C([0,T];BUC(\R^d;\C^N\times\C^N)))$ for all $|\alpha|=2m$,

$a_{\alpha} \in L^\infty(\O\times(0,T)\times \R^d;\C^{N}\times\C^N)$ for all $|\alpha|<2m$. \smallskip

\item
There is a constant $M_1\geq 0$ such that for all $t\in [0,T]$ and $\omega\in
\O$,
\[\sum_{|\alpha|=2m} \|a_{\alpha}(t,\omega,\cdot)\|_{\infty} \leq M_1.\]

\item There is a constant $M_2\ge 0$ and an angle $\vartheta\in [0,\frac12\pi)$
such that for all $t\in [0,T]$, $\omega\in \O$, $x\in \R^d$, and
$\xi\in \R^d$ with $|\xi|=1$ we have
$$\sigma(\mathcal{A}_{\pi}(t,\omega, x,\xi)) \subseteq \{z\in \C\setminus\{0\}: |\arg(z)|\leq\vartheta\}$$
and
\[ \|\mathcal{A}_{\pi}(t,\omega,x,\xi)^{-1}\|_{\calL(\C^N)}\leq M_2.\]
\end{enumerate}
\end{enumerate}
\let\theenumi\ALTERWERTA
\let\labelenumi\ALTERWERTB

\medskip
Let $A_q(t,\omega)$ denote the realization of
$\mathcal{A}(t,\omega,\cdot)$ in $L^q(\R^d;\C^N)$ with domain
\[\Dom(A(t,\omega)) = H^{2m,q}(\R^d;\C^N).\]
By \cite[Theorem 6.1]{DuSi}, applied pointwise on $\Omega$,
one has the following powerful result for the $H^\infty$-calculus
of $A$.

\begin{proposition}[\cite{DuSi}]
\label{prop:DuSi}
Let Hypothesis \ref{as:Acoef} be satisfied.
For all $q\in (1,\infty)$ and $\sigma\in (\vartheta, \frac12\pi)$
there exist constants $w\geq 0$ and $C\geq 1$,
depending only on $q$ $\sigma$, $\vartheta$, $M_1$, $M_2$,
such that for all $\omega\in \O$ and $t\in
[0,T]$ the operator $A_q(\omega,t)+w$ has a bounded $H^\infty(\Sigma_{\sigma})$-calculus on
$L^q(\R^d;\C^N)$ with boundedness constant at most $C$.
\end{proposition}

This result actually holds with $\vartheta\in [0,\pi)$, provided one extends
the definition of bounded $H^\infty$-calculi accordingly (replacing negative generators
of analytic semigroups by generals sectorial operators), but we shall not need it in
this generality.

\subsection{Hypotheses on the functions $f$, $f^0$, $b$, $b^0$,
and the initial value $u_0$}

\let\ALTERWERTA\theenumi
\let\ALTERWERTB\labelenumi
\def\theenumi{{\rm (Hf)}}
\def\labelenumi{(Hf)}
\begin{enumerate}
\item\label{as:eq1f} The function $f^0:[0,T]\times\O\times\R^d\to L^q(\R^d;\C^N)$ is
$\mathcal{P}\times \mathcal{B}_{\R^d}$-measurable and satisfies $f^0\in L^1(0,T;L^q(\R^d;\C^N))$ almost surely. The function $f:[0,T]\times\O\times\R^d\times H^{2m,q}(\R^d;\C^N)\to L^q(\R^d;\C^N)$ is
$\mathcal{P}\times \mathcal{B}_{\R^d}\times \mathcal{B}(H^{2m,q}(\R^d;\C^N))$-measurable.
There exist constants
$\alpha_f\in [0,1)$, $L_{f}\geq 0$, $L_{f,\alpha_f}\geq 0$, $C_f\geq 0$ such
that for all $u,v\in H^{2m,q}(\R^d;\C^N)$, $t\in [0,T]$, and $\omega\in \O$ one
has
\begin{align*}
\phantom{aaaa}
\|f(t,\omega,\cdot, u) & - f(t,\omega,\cdot, v)\|_{L^q(\R^d;\C^N)}\\ & \leq
 L_{f}
\|u-v\|_{H^{2m,q}(\R^d;\C^N)} + L_{f,\alpha_f} \|u-v\|_{H^{2m-\alpha_f,q}(\R^d;\C^N)}
\end{align*}
and
\begin{align*}
\phantom{aaaa}
\|f(t,\omega, u)\|_{L^q(\R^d;\C^N)} \leq C_{f}(1+ \|u\|_{H^{2m,q}(\R^d;\C^N)}).
\end{align*}
\end{enumerate}
\let\theenumi\ALTERWERTA
\let\labelenumi\ALTERWERTB

\let\ALTERWERTA\theenumi
\let\ALTERWERTB\labelenumi
\def\theenumi{{\rm (Hb)}}
\def\labelenumi{(Hb)}
\begin{enumerate}
\item\label{as:eq1b} The functions $b_{i}^0:[0,T]\times\O\times\R^d\to H^{m,q}(\R^d;\C^N)$ are
$\mathcal{P}\times \mathcal{B}_{\R^d}$-measurable and satisfy $b^0\in L^1(0,T;H^{m,q}(\R^d;\ell^2(\C^N)))$ almost surely.
The functions
$b_i:[0,T]\times\O\times\R^d\times H^{2m,q}(\R^d;\C^N)\to H^{m,q}(\R^d;\C^N)$ are
$\mathcal{P}\times \mathcal{B}_{\R^d}\times
\mathcal{B}(H^{2m,q}(\R^d;\C^N))$-measurable. There exist constants
$\a_b\in [0,1)$, $L_{b}\geq 0$, $L_{b,\a_b}\geq 0$ and $C_b$ such that
for all $u,v\in H^{2m,q}(\R^d;\C^N)$, $t\in [0,T]$, and $\omega\in \O$ one has
\begin{align*}
\phantom{aaaa}
\|b(t,\omega,\cdot, u) & - b(t,\omega,\cdot, v)\|_{H^{m,q}(\R^d;\ell^2(\C^N))}\\
& \leq  L_{b}
\|u-v\|_{H^{2m,q}(\R^d;\C^N)} + L_{b,\a_b} \|u-v\|_{H^{2m-\a_b,q}(\R^d;\C^N)}
\end{align*}
and
\begin{align*}
\phantom{aaaaa}
\|b(t,\omega, u)\|_{H^{m,q}(\R^d;\ell^2(\C^N))} \leq C_{b}(1+
\|u\|_{H^{2m,q}(\R^d;\C^N)}).
\end{align*}
\end{enumerate}
\let\theenumi\ALTERWERTA
\let\labelenumi\ALTERWERTB

\let\ALTERWERTA\theenumi
\let\ALTERWERTB\labelenumi
\def\theenumi{{\rm (H$u_0$)}}
\def\labelenumi{(H$u_0$)}
\begin{enumerate}
\item\label{as:equ0}
The initial value $u_0:\O\to L^q(\R^d;\C^N)$ is
$\F_0$-measurable.
\end{enumerate}
\let\theenumi\ALTERWERTA
\let\labelenumi\ALTERWERTB

\subsection{Main result}
We begin by defining the notion of a strong solutions to
the SPDE \eqref{eq:eq2m}. We fix exponents $p,q\in [2, \infty)$
and assume that \ref{as:Acoef}, \ref{as:eq1f}, \ref{as:eq1b}, \ref{as:equ0} are satisfied.
As in Section \ref{sec:4} it can be shown that a strong solution with paths
in $L^p(0,T;H^{2m,q}(\R^d;\C^N)))$ is also mild and weak solution (cf. Proposition
\ref{prop:strongmild} and the references given there).

\begin{definition}
A progressively measurable
process $u\in L^0(\O;L^p(0,T;\allowbreak H^{2m,q}(\R^d;\allowbreak \C^N)))$
is called a {\em strong solution} to \eqref{eq:eq2m} if for almost all
$(t,\omega)\in [0,T]\times \O$,
$$
\begin{aligned}
u(t,\cdot) +  \int_0^t \mathcal{A}(s,\cdot, D) u(s,\cdot)\, ds =  u_0(\cdot)   & +
\int_0^t f(s,\cdot,u(s,\cdot)) + f^0(s,\cdot)\,ds \\ & + \sum_{i\geq 1} \int_0^t
b_i(s,\cdot,u(s,\cdot)) + b_i^0(s,\cdot)\, d w_i(s).
\end{aligned}
$$
\end{definition}
The integral with respect to time is well defined as a
Bochner integral in the space
$L^q(\R^d;\C^N)$. By \eqref{eq:type2est} and the remark following it,
the stochastic integrals are well defined in the space $H^{m,q}(\R^d;\C^N)$.
Indeed, by \ref{as:eq1b} and the isomorphism \eqref{eq:BesselFubini}
one has
\begin{align*}
\big\|b(s,\cdot,u(s,\cdot))\big\|_{\g(\ell^2, H^{m,q}(\R^d;\C^N))} & \eqsim_q
\big\|b(s,\cdot,u(s,\cdot))\big\|_{H^{m,q}(\R^d;\ell^2(\C^N))}\\ & \leq C_{b}(1+
\|u(s,\cdot)\|_{H^{2m,q}(\R^d;\C^N)}).
\end{align*}
By the assumptions on $u$, the $L^2(0,T)$-norm of the right-hand side is
finite almost surely.

As a consequence of Theorem \ref{thm:SE2} one has the
following well-posedness result for the SPDE \eqref{eq:eq2m}.

\begin{theorem}\label{thm:eq2m}
Let $q\in [2, \infty)$ and $p\in (2, \infty)$, where $p=2$ is also allowed if
$q=2$. Assume \ref{as:Acoef}, \ref{as:eq1f},
\ref{as:eq1b}, \ref{as:equ0}, and suppose that
$f^0\in L^p_{\F}(\O;L^p(0,T;\allowbreak L^q(\R^d;\C^N)))$
and $b^0\in L^p_{\F}(\O;L^p(0,T; H^{m,q}(\R^d;\ell^2(\C^N))))$.
Provided $L_f$ and $L_b$ are small enough, the following assertions hold:

\begin{enumerate}[(i)]
\item If $u_0\in L^0_{\F_0}(\O;B^{2m(1-\frac1p)}_{q,p}(\R^d;\C^N))$, then the
problem \eqref{eq:eq2m} has a unique solution
$u\in L^0_{\F}(\O;L^p(0,T;H^{2m,q}(\R^d;\C^N)))$. Moreover, $u$ has a
version with trajectories in $C([0,T];B^{2m(1-\frac1p)}_{q,p}(\R^d;\C^N))$.

\item If $u_0\in L_{\F_0}^p(\Omega; B^{2m(1-\frac1p)}_{q,p}(\R^d;\C^N))$,
then the solution $u$ given by part {\rm (i)} satisfies
\begin{align*}
\|u\|_{L^p((0,T)\times\O;H^{2m,q}(\R^d;\C^N))} &\leq
C\big(1+\|u_0\|_{L^p(\O;B^{2m(1-\frac1p)}_{q,p}(\R^d;\C^N))}\big)\\
\phantom{aaaa} \|u\|_{L^p(\O;C([0,T];B^{2m(1-\frac1p)}_{q,p}(\R^d;\C^N))} & \leq
C\big(1+\|u_0\|_{L^p(\O;B^{2m(1-\frac1p)}_{q,p}(\R^d;\C^N))}\big),
\end{align*}
with constants $C$ independent of $u_0$.

\item For all $u_0, v_0\in
L_{\F_0}^p(\O;B^{2m(1-\frac1p)}_{q,p}(\R^d;\C^N))$,
the corresponding solutions $u,v$ satisfy
\begin{align*}
\|u-v\|_{L^p((0,T)\times\O;H^{2m,q}(\R^d;\C^N))} &\leq C
\|u_0-v_0\|_{L^p(\O;B^{2m(1-\frac1p)}_{q,p}(\R^d;\C^N))},\\
\phantom{aaaa} \|u-v\|_{L^p(\O;C([0,T];B^{2m(1-\frac1p)}_{q,p}(\R^d;\C^N)))}
& \leq C\|u_0-v_0\|_{L^p(\O;B^{2m(1-\frac1p)}_{q,p}(\R^d;\C^N))},
\end{align*}
with constants $C$ independent of $u_0$ and $v_0$.
\end{enumerate}
\end{theorem}
\begin{proof}
It suffices to check the conditions of Theorem \ref{thm:SE2} with $X_0
= L^q(\R^d;\C^N)$ and $X_1 = H^{2m,q}(\R^d;\C^N)$.
These spaces satisfy
Hypothesis \ref{as:XUMD}, Hypothesis \ref{as:A'} holds by Proposition \ref{prop:DuSi} and
the assumption that $\vartheta<\frac12\pi$, and Hypothesis \ref{as:initial_value} holds
by the assumption on $u_0$. The family
$\mathscr{J}$ is $R$-bounded from
$\calL(L^{p}_\mathscr{F}(\R_+\times\O ;\g(H,X_{0}))$ to
$L^{p}(\R_+\times\O ;X_0)$ by Theorem
\ref{thm:Rbddness}.

Recall from \cite[Theorems 2.4.2, 2.4.7 and 2.5.6]{Tri83} that
\begin{equation}\label{eq:interpspaces}
X_{\frac12} = H^{m,q}(\R^d;\C^N) \ \ \text{and} \ \  \Xp  =
B^{2m(1-\frac1p)}_{q,p}(\R^d;\C^N).
\end{equation}
Let $F:[0,T]\times\O\times X_1\to X_0$ be defined by $F(t,\omega,u) =
f(t,\omega,\cdot, u)$. The additional additive term can be defined in a similar way.
Then the equivalent version of \ref{as:LipschitzF} discussed in Remark \ref{rem:ass}
is satisfied with $\a_F =
\alpha_f$, $L_{F}' = L_{f}$, $\tilde L_{F}' = L_{f,\alpha_f}$ and $C_F = C_f$.
Let $H = \ell^2$ and let $B:[0,T]\times\O\times X_{1}\to \g(H,X_{\frac12})$ be
defined by $B(t,\omega,u) e_i =  b_i(t,\omega,\cdot, u)$. The additional additive
term can be defined in a similar way. Then the equivalent version of
\ref{as:LipschitzB} discussed in Remark \ref{rem:ass} is satisfied
with $\a_B = \a_b$, $L_{B}' = L_{b}$,
$\tilde L_{B}' = L_{b,\a_b}$ and $C_B = C_b$.

In this way, the equation \eqref{eq:eq2m} can be written as \eqref{SE'}, where
the unknown processes $u:[0,T]\times\O\times\R^d\to \C^N$ and
$U:[0,T]\times\O\to X_0$ are identified through $U(t,\omega)(x) =
u(t,\omega,x)$.
The result then follows from Theorem \ref{thm:SE2} and \eqref{eq:interpspaces}.
\end{proof}

\begin{remark}\label{rem:nZ}
Let $n\in \Z$. If $a_{\alpha}\in BUC^{|n|}(\R^d;\C^N\times\C^N)$ one can
transfer the result of Proposition \ref{prop:DuSi} to
the realization of $\mathcal{A}(t,\omega,\cdot)$ in $H^{n,q}(\R^d;\C^N)$ with
domain
\[\Dom(A_{n,q}(t,\omega)) = H^{n+2m,q}(\R^d;\C^N).\]
We refer to \cite[Lemma 5.2]{Kry} for details.
Using this fact, under suitably reformulated assumptions on $f$, $f^0$, $b$, $b^0$ and $u_0$
one can obtain a version of Theorem \ref{thm:eq2m} with an additional regularity
parameter $n\in \Z$. It is even possible to consider a real parameter $n$,
but in that case on needs additional smoothness on $a$ (see
\cite[Corollary 2.8.2]{Tri83}).
\end{remark}

\subsection{Discussion}
In this subsection we compare the above result Theorem \ref{thm:eq2m}
with available results in the literature.

The case $m=1$ and $N=1$ of Theorem \ref{thm:eq2m} has some overlap with
\cite[Theorem 5.1]{Kry} due to Krylov.
Theorem \ref{thm:eq2m} improves on
\cite[Theorem 5.1]{Kry} in various respects.

\begin{enumerate}
 \item[(i)] Our approach covers SPDEs governed by
$N$-dimensional systems of elliptic operators of order $2m$ for any $m\geq 1$.
\end{enumerate}

Even for $m=1$ and $N=1$, there are new features in our approach:

\begin{enumerate}
 \item[(ii)] In our setting,
the highest order coefficients $a_{\alpha}$ are only assumed to be bounded and
uniformly continuous in the space variable, whereas in \cite[Theorem 5.1]{Kry}
it is assumed that they are H\"older continuous in the space variable.
Our continuity assumptions can be further weakened to VMO
assumptions (cf.\ \cite{DuYa} for the second order case). Recently, in
\cite{Kim09} Krylov's $L^p$-approach has been extended to prove results for
continuous coefficients as well.
\item[(iii)]
In our approach, the parameters $p$ and $q$ can be chosen independently of each other.
In \cite[Theorem
5.1]{Kry}, only the case $p=q$ is considered, in \cite{Kry00} an extension to the case
$p\geq q\geq 2$ was obtained. We do not need such an assumption.
\end{enumerate}

Finally, the
regularity assumptions on the initial value in \cite[Theorem 5.1]{Kry} seem not
to be optimal.

On the other hand,
there are two striking features of Krylov's result that we could not
cover by our methods.
\begin{enumerate}
\item[(i)$'$] In \cite[Theorem 5.1]{Kry}, an additional
linear term satisfying a less restrictive smallness condition
can be allowed in the multiplicative part of the noise (see
\cite[Assumption 5.1]{Kry}).
\end{enumerate}
In our approach, we need a smallness condition on $L_f$ and $L_b$ and are
not able to take the linear part as mentioned above into account yet. There is a
possibility that the operator-theoretic approach of \cite{NVW2} works in such a
setting. We also refer to Subsection \ref{subsec:compH} for a discussion on the smallness
condition.
\begin{enumerate}
\item[(ii)$'$] In \cite[Theorem 5.1]{Kry}, the highest order coefficients
$a_{\alpha}$ with $|\alpha| = 2$ need only be measurable in time.
\end{enumerate}
Quite possibly, this cannot be achieved by an operator
theoretic approach. All well-posedness results for time-dependent problems
currently available in the literature
impose some continuity assumption in order to proceed by perturbation arguments.

With regard to (i), we mention that Mikulevicius and Rozovskii
\cite{MiRo} have extended Krylov's $L^p$-approach
to $N$-dimensional systems of second order equations. Apart from the fact that
our result covers operators of order $2m$, the differences
are of the same nature as those pointed out in (ii), (iii), and  (i)$'$, (ii)$'$.
A further difference is that  Mikulevicius and Rozovskii consider
equations in divergence form. Our results hold for systems of second
operators in divergence form as well, since, under mild regularity assumptions on the
coefficients, such operators also have a bounded $H^\infty$-calculus
(see \cite{DuMc, KuWe} and references therein).

\section{Second order parabolic SPDEs on bounded domains in $\R^d$}\label{sec:8}

We proceed with an application of Theorems \ref{thm:SE} and
\ref{thm:SE2} to a class of second order parabolic SPDEs on a
bounded domain $\mathcal{O}\subseteq \R^d$ with mixed Dirichlet and Neumann
boundary conditions. All results can be extended to
$N$-dimensional systems of operators of $2m$ for arbitrary $m\ge 1$,
assuming Lopatinskii-Shapiro boundary conditions (see \cite{DHP} for more on this).
The case $N=1$ and $m=1$ is chosen here in order to keep the technical details
at a reasonable level.

Let $\mathcal{O}\subseteq \R^d$ be a bounded domain with a $C^2$-boundary
$\partial \mathcal{O} = \Gamma_0 \cup \Gamma_1$ where $\Gamma_0$ and $\Gamma_1$
are disjoint and closed (one of them being possibly empty).
On $[0,T]\times \mathcal{O}$ we
consider the following stochastic partial differential equation with Dirichlet
boundary conditions on $\Gamma_0$ and Neumann boundary conditions on $\Gamma_1$:
\begin{equation}\label{eq:eqsecondorder}
\left\{\begin{aligned}
d u(t,x) +  \mathcal{A}(x, D) u(t,x)\, dt &= [f(t,x,u) + f^0(t,x)]\,dt
\\ & \qquad\qquad +
\sum_{i\geq 1} [b_i(t,x,u) + b_i^0(t,x)]\, d w_i(t),
\\ \mathcal{C}(x,D) u &=0,
\\ u(0,x) & = u_0(x).
\end{aligned}
\right.
\end{equation}
Here
\[\mathcal{A}(x, D) = \sum_{i,j=1}^d a_{ij}(x) D_i D_j + \sum_{i=1}^d a_i(x) D_i
+ a_0,\]
where $D_i$ denotes the $i$-th partial derivative, and
\[\mathcal{C}(x,D) = \sum_{i=1}^d c_{i}(x) D_i + c_0(x).\]

\subsection{Assumptions on the coefficients $a_{ij}$, $a_i$, $c_i$}
Essentially, the assumptions on $a_{ij}$ and $a_i$
correspond to a special case of an example in \cite{DDHPV} and \cite{KKW}.

\let\ALTERWERTA\theenumi
\let\ALTERWERTB\labelenumi
\def\theenumi{{\rm (Ha)}}
\def\labelenumi{(Ha)}
\begin{enumerate}
 \item\label{as:Acoefneumann} The coefficients $a_{ij}$, $a_i$, $c_i$ are
real-valued and satisfy:
\begin{enumerate}[(i)]
\item There is a constant $\rho\in (0,1]$ such that
\begin{align*}
a_{ij}& \in C^\rho(\overline{\mathcal{O}}) \ \  \text{for all $1\leq i,j\leq d$}.
\end{align*}
Furthermore,
\begin{align*}
   a_i & \in C(\overline{\mathcal{O}}) \ \ \text{for all $0\leq i\leq n$}, \\
   c_i & \in C^1(\overline{\mathcal{O}})  \ \ \text{for all $0\le i\le d$}.
\end{align*}

\item The matrices $(a_{ij}(x))$ are symmetric and there is a
constant $\kappa>0$ such that
for all $x\in \mathcal{O}$ and $\xi\in \R^d$ one has
\[\sum_{i,j=1}^d a_{ij}(x) \xi_i \xi_j \geq \kappa|\xi|^2.\]
\item
For all $x\in \Gamma_0$ we have $c_0(x) =1$ and $c_1(x) =c_2(x) = \ldots=c_d(x)= 0$.
There is a constant
$\kappa'>0$  such that for all $ x\in \Gamma_1$ we have
\[\sum_{i,j=1}^d c_{i}(x) n_i(x) \geq \kappa'.\]
\end{enumerate}
\end{enumerate}
\let\theenumi\ALTERWERTA
\let\labelenumi\ALTERWERTB

We denote by $A_q$ be the realization of $\mathcal{A}(\cdot)$ in
$L^q(\mathcal{O})$ with domain
\[\Dom(A(t,\omega)) = H^{2,q}_{\mathcal{C}}(\mathcal{O}):= \big\{u\in
H^{2,q}(\mathcal{O}): \mathcal{C}(x,D)u = 0\big\}.\]
One has the following result for
the $H^\infty$-calculus of $A_q$  (see \cite{DDHPV} and \cite{KKW}).

\begin{proposition}\label{prop:Hinftybdrneumann}
Assume that \ref{as:Acoefneumann} is satisfied.
For all $q\in (1, \infty)$ there exist
constants $w\geq 0$ and $\sigma\in [0,\frac12\pi)$ such that $A_q+w$ has a bounded
$H^\infty(\Sigma_{\sigma})$-calculus on $L^q(\mathcal{O})$.
\end{proposition}

\subsection{Hypotheses on the functions $f$, $f^0$, $b$, $b^0$,
and the initial value $u_0$}

\let\ALTERWERTA\theenumi
\let\ALTERWERTB\labelenumi
\def\theenumi{{\rm (Hf)}}
\def\labelenumi{(Hf)}
\begin{enumerate}
\item\label{as:eq1fneumann}
The function $f^0:[0,T]\times\O\times\OO\to L^q(\OO)$ is
$\mathcal{P}\times \mathcal{B}_{\OO}$-measurable and satisfies $f^0\in L^1(0,T;L^q(\OO))$ almost surely.
The function $f:[0,T]\times\O\times\mathcal{O}\times
H^{2,q}_{\mathcal{C}}(\mathcal{O})\to L^q(\mathcal{O})$ is $\mathcal{P}\times
\mathcal{B}_{\mathcal{O}}\times \mathcal{B}(H^{2,q}(\mathcal{O}))$-measurable
and there exist constants $\alpha_f\in [0,1)$, $L_{f}\geq 0$,
$L_{f,\alpha_f}\geq 0$, and $C_f\geq 0$ such that for all $u,v\in
H^{2,q}_{\mathcal{C}}(\mathcal{O})$, $t\in [0,T]$, and $\omega\in \O$ one has
\begin{align*}
\|f(t,\omega,\cdot, u) & - f(t,\omega,\cdot, v)\|_{L^q(\mathcal{O})}\\ & \leq
 L_{f}
\|u-v\|_{H^{2,q}(\mathcal{O})} + L_{f,\alpha_f} \|u-v\|_{H^{2-\alpha_f,q}(\mathcal{O})},
\end{align*}
and
\begin{align*}
\|f(t,\omega, u)\|_{L^q(\mathcal{O})} \leq C_{f}(1+\|u\|_{H^{2,q}(\mathcal{O})}).
\end{align*}
\end{enumerate}
\let\theenumi\ALTERWERTA
\let\labelenumi\ALTERWERTB

\let\ALTERWERTA\theenumi
\let\ALTERWERTB\labelenumi
\def\theenumi{{\rm (Hb)}}
\def\labelenumi{(Hb)}
\begin{enumerate}
\item\label{as:eq1bneumann}
The functions $b_{i}^0:[0,T]\times\O\times\OO\to L^q(\OO)$ are
$\mathcal{P}\times \mathcal{B}_{\OO}$-measurable and satisfy
$b^0\in L^1(0,T;H^{1,q}(\OO;\ell^2))$ almost surely.
The functions
$b_i:[0,T]\times\O\times\mathcal{O}\times H^{2,q}_{\mathcal{C}}(\OO)\to
H^{1,q}(\mathcal{O})$ are $\mathcal{P}\times \mathcal{B}_{\OO}\times
\mathcal{B}(H^{2,q}_{\mathcal{C}}(\OO))$-measurable and there exist constants
$\a_b\in [0,1)$, $L_{b,1}\geq 0$, $L_{b,\a_b}\geq 0$, and $C_b$ such that
for all $u,v\in H^{2,q}_{\mathcal{C}}(\OO)$, $t\in [0,T]$, and $\omega\in \O$
one has
\begin{align*}
\|b(t,\omega,\cdot, u) & - b(t,\omega,\cdot,
v)\|_{H^{1,q}_{\mathcal{C}}(\OO)}\\ & \leq
 L_{b} \|u-v\|_{H^{2,q}(\mathcal{O})} +
L_{b,\a_b}
\|u-v\|_{H^{2-\a_b,q}(\mathcal{O})}
\end{align*}
and
\begin{align*}
\|b(t,\omega, u)\|_{H^{1,q}_{\mathcal{C}}(\mathcal{O};\ell^2)} \leq C_{b}(1+
\|u\|_{H^{2,q}(\mathcal{O})}).
\end{align*}
\end{enumerate}
\let\theenumi\ALTERWERTA
\let\labelenumi\ALTERWERTB

\let\ALTERWERTA\theenumi
\let\ALTERWERTB\labelenumi
\def\theenumi{{\rm (H$u_0$)}}
\def\labelenumi{(H$u_0$)}
\begin{enumerate}
\item\label{as:equ00}
The initial value $u_0:\O\to L^q(\mathcal{O})$ is
$\F_0$-measurable.
\end{enumerate}
\let\theenumi\ALTERWERTA
\let\labelenumi\ALTERWERTB

\subsection{Main result}

We let $p,q\in [2, \infty)$ and assume that
 \ref{as:Acoef}, \ref{as:eq1f} \ref{as:eq1b}, \ref{as:equ00} are satisfied.

\begin{definition}
A progressively measurable process
$u\in L^0(\O;L^p(0,T;H^{2,q}(\mathcal{O})))$
is called a solution to \eqref{eq:eqsecondorder} if, for almost all $(t,\omega)\in [0,T]\times \O$,
\begin{equation*}
\begin{aligned}
u(t,\cdot) +  \int_0^t \mathcal{A}(\cdot, D) u(s,\cdot)\, ds =  u_0(\cdot)
& + \int_0^t f(s,\cdot,u(s,\cdot)) + f^0(s,\cdot)\,ds
\\ &  + \sum_{i\geq 1} \int_0^t
b_i(s,\cdot,u(s,\cdot)) + b^0_i(s,\cdot) \, d w_i(s).
\end{aligned}
\end{equation*}
\end{definition}
Arguing as in the previous section, we see that
the integral with respect to time is well defined as a Bochner integral in
the space $L^q(\mathcal{O})$ and the stochastic integrals are well defined
in $H^{1,q}_{\mathcal{C}}(\mathcal{O})$.

Following \cite{Amnonhom}, we define the following Besov and Bessel potential
spaces with boundary conditions. For $p\in (1, \infty)$ and $q\in (1, \infty)$,
and $S^s_q \in \{B_{q,p}^s, H^s_q\}$ let
\[S^{s}_{q,\mathcal{C}}(\mathcal{O}) = \{u\in S^{s}_{q}(\mathcal{O}):
\mathcal{C} u =0\} \ \ \ \text{$1+\tfrac1q<s\leq 2$}.\]
For $p\in (1, \infty)$ and $q\in (1, \infty)$, and $S^s_q \in \{B_{q,p}^s,
H^s_q\}$ let
\[S^{s}_{q,\mathcal{C}}(\mathcal{O}) = \{u\in S^{s}_{q}(\mathcal{O}):
\text{Tr}(u) =0 \ \text{on} \ \Gamma_0\} \ \ \ \text{$\tfrac1q<s<1+\frac1q$}. \]

Below we use the following well-known result:
\begin{equation*}
X_{\frac12} = H^{1,q}_{\mathcal{C}}(\mathcal{O}) \ \ \text{and} \ \  \Xp  =
B^{2-\frac2p}_{q,p,\mathcal{C}}(\mathcal{O}),
\end{equation*}
Indeed, since $A$ has a bounded $H^\infty$-calculus of angle $<\frac12\pi$,
it has bounded imaginary powers and therefore, by \cite[Theorem 1.15.3]{Tr1},
$X_{\frac12}  = [X_0, X_1]_{\frac12}$ with equivalent norms.
Now by Theorem 5.2 and Remark 5.3 (c) in \cite{Amnonhom} one obtains
$X_{\frac12} = H^{1,q}_{\mathcal{C}}(\mathcal{O})$  with equivalent norms.
Similarly, if ${2-\frac2p} \notin
\{\frac1q, 1+\frac1q\}$ then
\[\Xp  = (X_0, X_1)_{1-\frac1p,p} =
B^{{2-\frac2p}}_{q,p,\mathcal{C}}(\mathcal{O}).\]

Note that in the case that $\Gamma_0 = \emptyset$, one has
$H^{1,q}_{\mathcal{C}}(\mathcal{O}) = H^{1,q}(\mathcal{O})$  for all $q\in (1,
\infty)$  with equivalent norms.

As a consequence of Theorem \ref{thm:SE}
we obtain the following well-posedness result for the SPDE \eqref{eq:eq2m}.

\begin{theorem}\label{thm:eq2mneumann}
Let $q\in [2, \infty)$ and $p\in (2, \infty)$, where $p=2$ is also allowed if
$q=2$. Assume that $\frac{2}{p}+\frac1q\neq 1$.
Assume that \ref{as:Acoefneumann}, \ref{as:eq1fneumann},
\ref{as:eq1bneumann}, \ref{as:equ00} are satisfied and suppose that
$f^0\in L_{\F}^p(\O;L^p(0,T;L^q(\OO)))$ and $b^0\in L_{\F}^p(\O;\allowbreak L^p(0,T;\allowbreak H^{1,q}(\OO;\ell^2)))$.
Provided $L_f$ and $L_b$ are small enough, the following assertions
hold:

\begin{enumerate}[(i)]
\item If $u_0\in L_{\F_0}^0(\O; B^{{2-\frac2p}}_{q,p, \mathcal{C}}(\mathcal{O}))$,
then the problem \eqref{eq:eq2m} has a unique solution $u\in L^0_{\F}(\O;\allowbreak L^p(0,T;
\allowbreak H^{2,q}_{\mathcal{C}}(\mathcal{O})))$. Moreover, $u$
has a version with trajectories in the space $C([0,T];\allowbreak B^{{2-\frac2p}}_{q,p,
\mathcal{C}}(\mathcal{O}))$.
\item If $u_0\in L_{\F_0}^p(\Omega; B^{{2-\frac2p}}_{q,p,\mathcal{C}}(\mathcal{O}))$, then
the solution $u$ given by part {\rm (i)}
satisfies
\begin{align*}
\|u\|_{L^p((0,T)\times\O;H^{2,q}_{\mathcal{C}}(\mathcal{O}))} &\leq
C(1+\|u_0\|_{L^p(\O;B^{{2-\frac2p}}_{q,p,\mathcal{C}}(\mathcal{O}))})\\
\phantom{aaaa} \|u\|_{L^p(\O;C([0,T];B^{{2-\frac2p}}_{q,p,
\mathcal{C}}(\mathcal{O}))} & \leq
C(1+\|u_0\|_{L^p(\O;B^{{2-\frac2p}}_{q,p, \mathcal{C}}(\mathcal{O}))}),
\end{align*}
with constants $C$ independent of $u_0$.

\item For all $u_0, v_0\in L_{\F_0}^p(\O;B^{{2-\frac2p}}_{q,p,
\mathcal{C}}(\mathcal{O}))$,
the corresponding solutions $u,v$ satisfy
\begin{align*}
\|u-v\|_{L^p((0,T)\times\O;H^{2,q}_{\mathcal{C}}(\mathcal{O}))} &\leq C
\|u_0-v_0\|_{L^p(\O;B^{{2-\frac2p}}_{q,p, \mathcal{C}}(\mathcal{O}))},\\
\phantom{aaaa} \|u-v\|_{L^p(\O;C([0,T];B^{{2-\frac2p}}_{q,p,
\mathcal{C}}(\mathcal{O}))}
& \leq C\|u_0-v_0\|_{L^p(\O;B^{{2-\frac2p}}_{q,p, \mathcal{C}}(\mathcal{O}))},
\end{align*}
with constants $C$ independent of $u_0$ and $v_0$.
\end{enumerate}
\end{theorem}

\begin{proof}
We check the conditions of Theorem \ref{thm:SE2} with $X_0
= L^q(\mathcal{O})$ and $X_1 = H^{2,q}_{\mathcal{C}}(\mathcal{O})$.

As in the proof of Theorem \ref{thm:eq2m},
the verification of the  Hypotheses  \ref{as:XUMD},
\ref{as:A}, \ref{as:initial_value}, as well as
the $R$-boundedness of $\mathscr{J}$ is immediate.

Let $F:[0,T]\times\O\times X_1\to X_0$ be defined by $F(t,\omega,u) =
f(t,\omega,\cdot, u)$. The additional term can be defined in a similar way.
Then the equivalent version of \ref{as:LipschitzF} discussed in Remark \ref{rem:ass} is satisfied with $\a_F =
\alpha_f$, $L_{F}' = L_{f}$, $\tilde L_{F}' = L_{f,\alpha_f}$ and $C_F = C_f$.
Let $H = \ell^2$ and let $B:[0,T]\times\O\times X_{1}\to \g(H,X_{\frac12})$ be
defined by $B(t,\omega,u) e_i =  b_i(t,\omega,\cdot, u)$. The additional term can be defined in a similar way. Then
\ref{as:LipschitzB} (see Remark \ref{rem:ass}) is satisfied with $\a_B = \a_b$, $L_{B}' = L_{b}$,
$\tilde L_{B}' = L_{b,\a_b}$ and $C_B = C_b$.

In this way, the equation \eqref{eq:eq2m} can be written as \eqref{SE'}, where
the unknown processes $u:[0,T]\times\O\times\mathcal{O}\to \R$ and
$U:[0,T]\times\O\to X_0$ are identified through $U(t,\omega)(x) =
u(t,\omega,x)$.
The result now follows from Theorem \ref{thm:SE2} and \eqref{eq:interpspaces}
and the assumptions on $p$ and $q$.
\end{proof}

\begin{remark}
Under additional continuity assumptions on the
coefficients $a_{ij}$, the same methods one can be used to
handle the case where $\mathcal{A}$ depends on time and $\O$.
\end{remark}

\subsection{Discussion} \
In case of Dirichlet boundary conditions, related results
for weighted half-spaces and bounded domains with weights
have been obtained by Kim and Krylov (see
\cite{Kim09} and references therein) using Krylov's $L^p$-approach.
The weighted approach
started with the $L^2$-theory of Krylov \cite{Kry94b}.
The advantage of using weights is that no additional compatibility conditions on the
noise are required in this case, whereas in the unweighted case such
conditions seem to be unavoidable (see \cite{Fla90}).
To see the point, note that Theorem \ref{thm:eq2mneumann} does not cover the simple problem
$$
\left\{
\begin{aligned}
du(t,x) & = \tfrac12 \Delta u(t,x)\,dt + dw(t), && t\in [0,T], \ x\in (0,1), \\
 u(t,0) = u(t,1) & = 0,                             && t\in [0,T], \\
 u(0,x) & = 0,                             && x\in (0,1),
\end{aligned}
\right.$$
in, say, $E = L^q(0,1)$ with $q\in (2,\infty)$, where $w$ is a real-valued Brownian motion.
Now  the constant function $g= \one $ in the noise term $dw(t) = \one\, dw(t)$ does not belong to
$\Dom((-\Delta)^\frac12) = H_0^{1,q}(0,1)$ due to the Dirichlet boundary
conditions. We refer to \cite[Section 4]{KryOverview}
for a further discussion of this example.

It seems likely that, under suitable regularity assumptions on the
coefficients, the $L^p$-realisations of the operators $\mathcal{A}$
on weighted domains should have a bounded $H^\infty$-calculus.
If true,
the weighted domain case could be treated by our methods as well.
Maximal $L^p$-regularity results
for elliptic operators on $L^p(\mathcal{O},w)$ for
open domains $\mathcal{O}\subseteq \R^d$ and Muckenhaupt weights $w$
were proved in \cite{HHH}.
Further evidence  is provided by the fact
(see \cite{DevNWe}) that if
the linear Cauchy problem with additive noise has maximal regularity
for both $A$ and its adjoint, then $A$ necessarily has a bounded
$H^\infty$-calculus.

\section{The stochastic Navier-Stokes equation}\label{sec:NS}

Let $d\ge 2$ be a fixed integer and suppose that $\OO$ is a smooth bounded open domain in $\R^d$. Let $H$ be a Hilbert space (for instance $\ell^2$ or $L^2(\OO)$).
Let $q\in (1,\infty)$ be fixed.
We are interested in local existence of strong solutions in $(H^{1,q}(\OO))^d$ of
the Navier-Stokes equation
\begin{equation}\label{eq:NS}
\left\{
\begin{aligned}
\frac{\partial u}{\partial t} & = \Delta u - (u\cdot \nabla) u + f^0 - \nabla p  + (g(u,\nabla u) + g^0) \dot W_H, \\
\div u(t,\cdot) & = 0, \quad t>0,\\
u(t,x) & = 0, \quad t>0, \ x\in \partial \OO,\\
u(0,\cdot) & = u_0.
\end{aligned}
\right.
\end{equation}
Note that we allow $g$ to depend on both $u$ and $\nabla u$.
As is well known (see, for instance,
\cite{BrCaFl, MiRo04}) such dependencies arise in the
modelling of the onset of turbulence.

The function $u_0:\OO\to \R^d$ is the initial velocity field,
$W_H$ is a cylindrical Brownian motion in $H$,
and $u$ and $p$ represent the velocity field and the pressure of the fluid, respectively.
We assume that $f^0$ and $g^0$
are strongly measurable and adapted and belong to $L^1(0,T;H^{-1,q}(\OO))$
and $L^2(0,T;L^q(\OO;H))$ almost surely, respectively. The function
$g$ is interpreted as a strongly measurable mapping
$$g: (H^{1,q}(\OO))^d\to(L^q(\OO;H))^d,$$
and we assume that
for and all $x,y\in (H^{1,q}(\OO))^d$
we have
\begin{equation}\label{eq:Lip-g}
\|g(u) - g(v)\|_{(L^q(\OO;H))^d}\leq L_g \|u-v\|_{(H^{1,q}(\OO))^d} + \tilde{L}_{g} \|u-v\|_{(L^q(\OO))^d}.
\end{equation}

It is well known (see \cite{FujMor} and \cite[Section 9]{KKW})
that we have the direct sum decomposition
$$(L^q(\mathcal{O}))^d = \X \oplus {\mathbb{G}^q},
$$ where $\X$ is the closure in $(L^q(\mathcal{O}))^d$ of the
set $\{u\in (C_{\rm c}^\infty(\mathcal{O}))^d: \ \nabla \cdot u = 0\}$ and
${\mathbb{G}^q} = \{\nabla p: \ p\in H^{1,q}(\mathcal{O})\}$.
We denote by $P$ the Helmholtz projection
from $(L^q(\mathcal{O}))^d$ onto $\X$ along this decomposition.
The negative
{\em Stokes operator} is the linear operator $(A,\Dom(A))$ defined by
\begin{align*}
\Dom(A)  & = \X \cap \Dom(\Delta_{\rm Dir}), \\
 Av   & = -P(\Delta u), \quad  u\in \Dom(A),
\end{align*}
where $\Dom(\Delta_{\rm Dir})$ is the domain of the Dirichlet Laplacian in $(L^q(\OO))^d$, which for $C^{2}$-domains equals
\[\Dom(\Delta_{\rm Dir}) =\{u\in (H^{2,q}(\OO))^d: \ u = 0 \hbox{ on } \partial \OO \}.\]
The operator $A$ is boundedly invertible (see \cite{Catt} and \cite[page 797]{KKW}), $-
A$ generates a bounded
analytic $C_0$-semigroup in $\X$, and it was shown in \cite{NoSa} (for $C^3$ domains)
and \cite[Theorem 9.17]{KKW} (for $C^{1,1}$ domains) that $A$ has a
bounded $H^\infty$-calculus on $\X$:

\begin{proposition}\label{prop:hinftystokes}
For all $q\in (1, \infty)$
the negative Stokes operator $A$ has a bounded $H^\infty(\Sigma_{\sigma})$-calculus
of angle $0<\sigma<\frac12\pi$ on $\X$.
\end{proposition}

It is well known (see \cite{Sohr} for the details) that, by
applying the Helmholtz projection $P$ to $u$,
the Navier-Stokes equation \eqref{eq:NS} can be reformulated as an abstract
stochastic evolution on $$ X_0 := \X_{-\frac{1}{2}},
$$ where the space on the right-hand side is defined as the completion of
$\X$ with respect to the norm
$$ \n x\n_{X_0} := \n A^{-\frac12}x\n_{\X}.$$
In particular, as a Banach space, $X_0$ is isomorphic to a closed subspace
of $L^q(\OO)$.

The bounded invertibility of $A$ implies that the identity
operator on $\X$ extends to a
continuous embedding $\X \embed X_0 $.
Furthermore we set
$$ X_1 := \Dom(A^\frac12).$$

For $s\in (0, 1]$  and $s-\frac1q>0$ let $H^{s,q}_0(\OO)$ and $B^{s}_{q,p,0}(\OO)$
denote the closed subspaces of $H^{s,q}(\OO)$ and $B^{s}_{q,p}(\OO)$ with zero trace. If $s-\frac1q<0$ we let
$H^{s,q}_{0}(\OO) =  H^{s,q}(\OO)$ and
$B^{s}_{q,p,0}(\OO) = B^{s}_{q,p}(\OO)$. Furthermore let $H^{-1,q}(\OO)$
be the dual of $H^{1,q'}_0(\OO)$ with $\frac1q + \frac1{q'} =1$.

The following lemma is well known.
\begin{lemma}\label{lem:interpostokes}
For every $\alpha\in [\frac12, 1]$ and $p,q\in (1, \infty)$ with $2\alpha-1-\frac1q \neq 0$ one has
\begin{equation}\label{eq:Xalphastokescomplex}
X_{\alpha} = \Dom(A^{\alpha-\frac12}) = \X\cap (H^{2\alpha-1,q}_{0}(\OO))^d,
\end{equation}
\begin{equation}\label{eq:Xalphastokesreal}
X_{\alpha,p} = \X\cap (B^{2\alpha-1}_{q,p,{0}}(\OO))^d.
\end{equation}
Moreover,
$P$ induces a bounded linear operator
\begin{align*}
 & P:(H^{-1,q}(\OO))^d\to X_0.
\end{align*}
\end{lemma}
\begin{proof}
To prove \eqref{eq:Xalphastokescomplex} note that
\begin{equation}\label{eq:Xalphastokes}
X_{\alpha} = \Dom(A^{\alpha-\frac12}) = [\X, \Dom(A)]_{\alpha-\frac12},
\end{equation}
where we used \cite[Theorem V.1.5.4]{Am}, Proposition \ref{prop:hinftystokes}
and \eqref{eq:fracpow}. The second identity in \eqref{eq:Xalphastokescomplex}
follows from \cite{Gigafract}, \cite[Theorem 9.17]{KKW} and
\cite[Theorem 1.17.1.1]{Tr1}. By a similar reasoning one obtains
\eqref{eq:Xalphastokesreal}.

The final assertion follows from a similar argument as in \cite[Proposition 9.14]{KKW} (see also \cite[Proposition 3.1]{Ku10}).
\end{proof}

Define $F:X_{\theta+\frac12}\times X_{\theta+\frac12}\to X_0$
by
$$ F(u,v) = -P((u\cdot \nabla)v)$$
and write $F(u):= F(u,u).$
We will check that these mappings are well defined for $\theta\geq \frac{d}{4q}$.
Indeed, by \cite{GigMiy}, for these $\theta$ one has
\[ \| A^{-\frac{1}{2}} F(u,v) \|_{L^q(\OO)} \le C \|A^{\theta} u \|_{L^q(\OO)} \| A^{\theta} v\|_{L^q(\OO)}, \ \ \ u,v\in \Dom(A^{\theta}).\]
This can be reformulated as
\[\|F(u,v)\|_{X_{0}}\leq C\|u\|_{X_{\theta+\frac12}}\|v\|_{X_{\theta+\frac12}}, \ \ \ u,v\in X_{\theta+\frac12},\]
from which the well-definedness follows.
Moreover, one immediately obtains the following local Lipschitz estimate (see \cite{BrzPes99})
\[  \| A^{-\frac{1}{2}} (F(u) - F(v)) \|_{L^q(\OO)} \le C (\| A^{\theta} u\|_{L^q(\OO)} + \| A^{\theta} v\|_{L^q(\OO)}) \| A^{\theta} u - A^{\theta} v\|_{L^q(\OO)},\]
which can be reformulated as
\[\|F(u) - F(v)\|_{X_{0}}\leq C (\|u\|_{X_{\theta+\frac12}}+ \|v\|_{X_{\theta+\frac12}})
\|u-v\|_{X_{\theta+\frac12}}.\]
In particular, if $0\le \theta<\frac{1}{2}-\frac1p$
and $p\in [2, \infty)$, then $1-\frac1p>\theta+\frac12$ and therefore
$F: \Xp\times\Xp\to X_0$ is locally Lipschitz continuous.

Next define $B:X_{1}\to \g(H,X_{\frac{1}{2}})$ by
$$ B(u) = P(g(u)).$$
This is well defined, because $g$ maps $X_1 = (H^{1,q}(\OO))^d$ into
\[(L^q(\OO;H))^d
= (\g(H,L^q(\OO))^d =  \g(H,(L^q(\OO))^d),\]
and the Helmholtz projection extends to a bounded projection
$P:\g(H,(L^q(\OO))^d)\to \g(H,\X) = \g(H,X_{\frac{1}{2}})$ in a canonical way. Here we used \eqref{eq:gammaFubini0}
and \eqref{eq:Xalphastokes}.

Now we can reformulate \eqref{eq:NS} as an abstract stochastic evolution
equation in $X_0$ of the form
\begin{equation}\label{eq:NS-abstract}
\left\{
\begin{aligned}
dU(t) + \mathscr{A} U(t)\,dt  & = [F(U(t)) + f(t)] \,dt + [B(U(t)) + b(t)]\,dW_H(t), \\
              U(0) & = u_0,
\end{aligned}
\right.
\end{equation}
where $\mathscr{A} = A_{-\frac12}$, $f = P f^0$ and $b = P g^0$.

\begin{theorem}
Let $d\ge 2$, and let $p>2$ and $q\ge 2$ satisfy
$\frac{d}{2q} < 1-\frac2p$.
Let $u_0: \Omega\to \X\cap B^{1-\frac2p}_{q,p,0}(\OO))^d$
be strongly $\F_0$-measurable. Let $f^0\in L^0_{\F}(\O;L^p(0,T;(H^{-1,q}(\OO))^d))$ and $g^0\in L^0_{\F}(\O;L^p(0,T;L^q(\OO;H)))$.
If the Lipschitz constant $L_g$ in \eqref{eq:Lip-g} is small enough, then
the problem \eqref{eq:NS-abstract} admits a unique
maximal local mild solution on $[0,T]$ with values in $(H_0^{1,q}(\OO))^d$. Moreover,
this solution has a modification with continuous trajectories in
$(B^{1-\frac2p}_{q,p,0}(\OO))^d$.
\end{theorem}
\begin{proof}
The operator family
$\mathscr{J}$ is $R$-bounded. Furthermore, by Proposition \ref{prop:hinftystokes}
$A$ has a bounded $H^\infty$-calculus on $\X = X_{\frac12}$ (the equality
of these spaces follows from \eqref{eq:Xalphastokescomplex}) of angle $<\frac12\pi$.
Therefore, $\mathscr{A} = A_{-\frac12}$ has a bounded $H^\infty$-calculus on $X_0$.

By Lemma \ref{lem:interpostokes}
(and noting that $1-\frac2p>\frac{d}{2q}\ge \frac1q$ to justify the
boundary conditions),
one has $X_0 = \X\cap (H^{-1,q}(\OO))^d$ $ X_1 = \X \cap (H^{1,q}_{0}(\OO))^d$ and
\[(X_0, X_1)_{1-\frac1p, p} = \X \cap (B^{1-\frac2p}_{q,p,0}(\OO))^d,
\ \ \ \ (X_0, X_1)_{\frac12} = \X.\]
By Lemma \ref{lem:interpostokes}, $u_0\in (X_0, X_1)_{1-\frac1p, p}$ almost surely.

For any $\theta\in [\frac{d}{4q},\frac12-\frac1p)$,
we can apply Theorem \ref{thm:SElocal} with $F^{(1)} = 0$, $F^{(2)}=F$,
 $B^{(1)} = B$, and $B^{(2)} = 0$ (and combine \eqref{eq:Lip-g} with Remark \ref{rem:ass} to check the
assumptions concerning $B^{(1)}$)
to obtain a unique maximal local mild solution $U$ which
satisfies the assertions of Theorem \ref{thm:SElocal}.
\end{proof}

\begin{remark}
The above result is merely a proof-of-principle and can
be extended into various directions. For instance,
more general ranges of the parameters can be considered
as in \cite{BrzPes99, GigMiy}; different regularity assumptions on the coefficients are possible,
and different regularity of the solutions will result. Furthermore,
we expect global existence in dimension $d=2$.
Using the results of \cite{Ku08, Ku10}, we believe that it should be possible
to adapt the above techniques to study maximal regularity for
the Navier--Stokes equation on $\R^d$ (see also the discussion below).
Along similar lines, it should be possible to use
the results of \cite{Ku08, Ku10, NoSa} to
study maximal regularity in the case of exterior
domains in $\R^d$. We plan to address such issues in a
forthcoming paper.
\end{remark}

\subsection{Discussion}

The existence of $H^{1,q}(\OO)$-solutions for the stochastic Navier-Stokes equation
in dimension $d=2$
was established, under a trace class assumption on the noise replacing our assumption
on $g$,
by Brze\'zniak and Peszat \cite{BrzPes99}. In their framework,
$g$ is a $C^1$-function on $\R^d$ with locally Lipschitz continuous derivatives;
it is then shown that $g$ induces a locally Lipschitz continuous mapping
$G$ from $X_\eta$ to
$\gamma(H,X_\frac12)$ for suitable exponents $\eta>1$. However, $G$
is not defined on $X_1$ and therefore
$g$ cannot be allowed to depend on both $u$ and $\nabla u$.

Under the same assumptions on $g$ as ours, existence of a local
strong $H^{1,q}(\R^d)$-solution for dimensions $d\ge 2$ has been
shown by Mikulevicius and Rozovskii \cite{MiRo04}.
Existence and uniqueness of local strong
$H^{1,2}$-solutions in bounded domains
was obtained by Mikulevicius \cite{Mik09}. In both papers,
global existence for $d=2$ is established as well.

\section*{Acknowledgments} We thank the anonymous referees for their detailed and helpful comments.

\def\polhk#1{\setbox0=\hbox{#1}{\ooalign{\hidewidth
  \lower1.5ex\hbox{`}\hidewidth\crcr\unhbox0}}} \def\cprime{$'$}

\end{document}